\documentclass[reqno]{amsart}
\usepackage[a4paper,centering,vscale=0.75,hscale=0.75]{geometry}
\usepackage[T1]{fontenc}
\usepackage[latin1]{inputenc}
\usepackage[english]{babel}
\usepackage[babel]{csquotes}

\usepackage[colorlinks=true, pdfstartview=FitV, linkcolor=blue, citecolor=blue, urlcolor=blue]{hyperref}

\usepackage{cite}
\usepackage{amssymb}
\usepackage{amsmath}
\usepackage{amsthm}
\usepackage{latexsym}
\usepackage{graphicx}
\usepackage{mathrsfs}
\usepackage{mathrsfs}
\usepackage{bbm}
\usepackage{verbatim}

\numberwithin{equation}{section}

\usepackage[usenames,dvipsnames]{color}

\newtheorem{thm}{Theorem}[section]

\newtheorem{lem}[thm]{Lemma}

\newtheorem{defin}[thm]{Definition}
\newtheorem{remark}[thm]{Remark}

\def\enne{\mathbb{N}}

\def\erre{\mathbb{R}}

\def\P{\mathbb{P}}

\def\E{\mathop{{}\mathbb{E}}}
\def\mN{\mathcal N}
\def\cL{\mathscr{L}}
\def\cF{\mathscr{F}}

\def\eps{\varepsilon}
\def\cP{\mathscr{P}}

\def\OO{\mathcal{O}}
\renewcommand{\div}{\operatorname{div}}
\renewcommand{\d}{{\mathrm d}}
\def\bu{{\mathbf u}}
\def\bv{{\mathbf v}}
\def\bh{{\mathbf h}}
\def\mU{\mathcal U}

\def\beq{\begin{equation}}
\def\eeq{\end{equation}}

\def\to{\rightarrow}
\def\wto{\rightharpoonup}
\def\wstarto{\stackrel{*}{\rightharpoonup}}
\def\embed{\hookrightarrow}

\def\norm #1{\left\|#1\right\|}

\def\sp #1#2{\left<#1,#2\right>}
\newcommand\ip\sp

\begin{document}
\title[Stochastic convective Cahn-Hilliard equation]
{Analysis and optimal velocity control
of\\ a stochastic convective Cahn-Hilliard equation}

\author{Luca Scarpa}
\address{Faculty of Mathematics, University of Vienna
Oskar-Morgenstern-Platz 1, 1090 Vienna, Austria.}
\email{luca.scarpa@univie.ac.at}
\urladdr{http://www.mat.univie.ac.at/$\sim$scarpa}

\subjclass[2010]{35K25, 35R60, 60H15, 80A22}
\keywords{Stochastic Cahn-Hilliard equation; convection; well-posedness;
optimal velocity control; optimality conditions.}   

\begin{abstract}
  A Cahn-Hilliard equation with stochastic multiplicative noise
  and a random convection term is considered.
  The model describes isothermal phase--separation occurring in a moving fluid,
  and accounts for the randomness appearing at the microscopic level
  both in the phase--separation itself and
  in the flow--inducing process.
  The call for a random component in the convection term 
  stems naturally from applications, as 
  the fluid's stirring procedure 
  is usually caused by mechanical or magnetic devices.
  Well--posedness of the state system is addressed and
  optimisation of a standard tracking type cost
  with respect to the velocity control is then studied.
  Existence of optimal controls is proved and 
  the G\^ateaux--Fr\'echet differentiability 
  of the control-to-state map is shown.
  Lastly, the corresponding adjoint backward problem is analysed, and
  first-order necessary conditions for optimality are derived
  in terms of a variational inequality involving the intrinsic adjoint variables.
\end{abstract}

\maketitle

%%%%%%%%%%%%%%%%%%%%%%%%%%%%%%%%%%%%%%%%%%%%%%%%

\section{Introduction}
\setcounter{equation}{0}
\label{sec:intro}
The aim of this paper is to analyse 
the stochastic Cahn-Hilliard equation with convection
\begin{align}
  \label{eq1}
  \d \varphi - \Delta\mu\,\d t + \bu\cdot\nabla\varphi\,\d t =
  B(\varphi)\,\d W \qquad&\text{in } (0,T)\times\OO=:Q\,,\\
  \label{eq2}
  \mu = -\Delta\varphi + \Psi'(\varphi) \qquad&\text{in } (0,T)\times\OO\,,\\
  \label{eq3}
  {\bf n}\cdot\nabla\varphi = {\bf n}\cdot\nabla\mu= 0 
  \qquad&\text{in } (0,T)\times\partial\OO\,,\\
  \label{eq4}
  \varphi(0)=\varphi_0 \qquad&\text{in } \OO\,,
\end{align}
where $\OO$ is a smooth bounded domain in $\erre^d$, $d=2,3$,
$T>0$ is a fixed final time, and $\bf n$ denotes the 
normal outward unit vector on $\partial\OO$.
The system \eqref{eq1}-\eqref{eq4} models 
isothermal phase--separation occurring in a moving fluid 
occupying the space region $\OO$ during the time interval $[0,T]$.
The order parameter, or phase--variable, $\varphi$ 
represents the relative concentration 
between the pure phases,
the variable $\mu$ represents the chemical potential of the system,
and the nonlinearity $\Psi:\erre\to\erre$ is a double-well potential
with two global minima. The term $\bu$
is an external random velocity field acting on the system,
modelling possible stirring and mixing processes of the fluid
which may affect phase--separation itself.
The stochastic forcing describing the 
thermal fluctuations affecting phase--separation is modelled 
by means of a cylindrical Wiener process $W$
on a given probability space and a $W$-integrable coefficient $B$,
possibly depending on the phase variable itself, which calibrates
the intensity of the noise.

The Cahn--Hilliard equation is a classical 
model employed in phase--separation, and has nowadays
numerous applications to physics, biology, and engineering. 
Its introduction 
dates back to the pioneering work by Cahn \& Hilliard \cite{cahn-hill}, 
where it was proposed, in the deterministic version,
to adequately describe spinodal decomposition in binary metallic alloys.
In the last decades the model has been extensively refined
in several directions. For example, 
the description of possible viscous behaviours has been 
originally presented in \cite{ell-st, ell-zhen, novick-cohen},
and then generalised in \cite{gurtin}.
The presence of a further evolution close to boundary
due to the interaction with the hard walls
has been accounted for by proposing 
several choices of dynamic boundary conditions, 
for which we refer to \cite{fish-spinod, kenz-spinod, gal-DBC2}.

The deterministic Cahn--Hilliard equation has been 
proven to be extremely effective in describing phase--separation phenomena.
Nevertheless, it presents some drawbacks. Indeed, 
the phase--separation process inevitably presents some disruptions,
acting at a microscopic level. These are due to
unpredictable movements at the atomistic level, which may be caused
for example by temperature oscillations, magnetic effects,
or configurational interactions. 
As such, the classical Cahn--Hilliard system is unable to capture 
the erratic nature of the separation process.
The most natural way to overcome this problem 
is to switch to a random setting instead, by introducing a
suitable noise term in the equation that could 
effectively describe the unpredictability of the phenomenon 
at a small scale. This was proposed by Cook in \cite{cook} for Wiener--type noises
and gave rise to the well--known Cahn--Hilliard--Cook
stochastic model for phase--separation.
The stochastic version of the model was then confirmed 
multiple times \cite{bin-SCH, pego-SCH} 
to be the only one that can genuinely describe phase--separation in alloys.
Since then, the random version of the equation has been 
increasingly studied, both in the physics literature 
\cite{rer-SCH, rer-SCH2, gsvg-SCH, lbh-SCH, mhb-SCH}
and in the direction of 
model validation and numerical simulations
\cite{bmw-SCH, bmw-SCH2, bmw-SCH3, haw-SCH, haw-SCH2, haw-SCH3, lee-CH}.

The classical  Cahn--Hilliard equation 
is the gradient flow associated to
the free energy functional 
\[
  \varphi\mapsto \frac12\int_\OO|\nabla\varphi|^2 + \int_\OO \Psi(\varphi)\,,
\]
with respect to the metric of $H^1(\OO)^*$. The gradient 
term penalises the oscillation of the order parameter, 
while the double--well potential models the tendency of each 
phase to concentrate. The form of the chemical potential
in \eqref{eq2} appears then 
naturally from the differentiation of the free energy.
Typical examples of $\Psi$
are given by 
\beq\label{psi_log}
  \Psi_{log}(r):=
  \frac{\theta}{2}\left((1+r)\ln(1+r) + (1-r)\ln(1-r)\right) 
  - \frac{\theta_0}{2}r^2\,, \quad r\in(-1,1)\,, \qquad 0<\theta<\theta_0\,,
\eeq
and
\beq\label{psi_pol}
  \Psi_{pol}(r):=\frac14(r^2-1)^2\,, \quad r\in\erre\,.
\eeq
Although \eqref{psi_log} is the most relevant choice 
in terms of thermodynamical consistency, its singular 
behaviour in $\pm1$ could be hard to tackle from the 
mathematical viewpoint, and in several 
models the polynomial approximation \eqref{psi_pol} is often employed.

The velocity field $\bu$ models the transport effects due to 
convection terms acting on the system. 
In our analysis, this will be a prescribed external 
forcing field which will play the role of velocity control
in a typical optimisation problem.
Optimisation involving phase-separating fluids 
where the velocity is the control arises naturally in applications.
For example, this is the case of 
block solidification of silicon crystals in photovoltaic applications: here,
the flow of the fluid acts as a control to optimise the 
distribution of certain impurities, at the atomistic level, in a process 
of solidification of silicon melt. For more details about the 
applications of optimal velocity control problem in phase--separating fluids
we refer to \cite{kud, roc-spr}.
In practice, the motion of the fluid can be 
achieved in several ways: as 
pointed out in \cite{col-gil-spr-conv, roc-spr}, 
the most common choices consist in
employing either mechanical stirring devices 
or ultrasound emitters directly into the container.
Another possibility is 
to prescribe a velocity on the fluid by means 
of magnetic fields: this is widely employed 
for example in the case of 
molten metals \cite{kud} or bulk semiconductor crystals.
Nevertheless, it is worthwhile noting that in all 
these scenarios, the velocity field
is usually obtained in an {\em indirect} way, meaning that 
the motion of the fluid is achieved only as a consequence 
of more {\em direct} controls, such as
mechanical devices or magnetic effects.
This being noticed, it is clear then that the the external prescription 
of a given velocity is strongly affected by microscopic noises,
which may be caused, depending on the type of 
motion--inducing devices, by configurational 
or electromagnetic disturbances occurring in the flow--creating process.
Also, the effective induction of the flow is strongly affected
by the imprecision of the above--mentioned devices.
From the modelling point of view, this strongly calls
for the introduction of a further source of randomness in the velocity 
field $\bu$ and for abandoning the classical deterministic setting of the problem.
Let us stress that the random component of the velocity field 
prescinds from the stochastic nature of the noise in equation \eqref{eq1}:
while the Wiener process $W$ models microscopic turbulences 
occurring in phase--separation, the random nature of $\bu$
takes into account the imprecision of the flow--inducing mechanisms.
For example, in typical situations $\bu$ would satisfy a further 
stochastic equation involving a further Wiener process, independent of $W$.
Clearly, this extra equation would specifically depend on the model in consideration:
here, in order to make the treatment as general and light as possible,
we only require $u$ to be a stochastic process.
Besides, the importance of allowing the control variable to be random 
is crucial when dealing with a controlled stochastic equation (see e.g.~\cite{yong-zhou}).
Indeed, bearing in mind the typical perspective of Monte Carlo simulations, 
restricting to deterministic controls would mean to choose {\em a priori}
a control which is independent of the possible outcomes of the evolution
according to the prescribed underlying probability space. By contrast,
stochastic controls ensure more freedom from the point of view of the controller,
as they allow to adapt the control to the random outcomes of the 
phenomenon itself.
With this in mind, in our analysis $\bu$ will be a prescribed stochastic process
satisfying some natural box--constraints, possibly taking into 
account the random imprecision of the velocity--inducing devices.
The model that we study presents then two main sources of randomness:
the first one is given by the Wiener noise in equation \eqref{eq1},
taking into account the microscopic turbulence affecting 
phase--separation, and the second one is the 
stochastic component of the convection term,
modelling the imprecision of the stirring procedure.
Hence, one can think the two random forcings as acting on 
two separate levels: a microscopic scale described by $W$,
and a different uncorrelated scale rendered by $\bu$.

The mathematical literature dealing with the Cahn-Hilliard equation 
is extremely developed.
In the deterministic case, attention has been widely devoted to 
the study of well-posedness, regularity, long--time behaviour of solutions,
and asymptotics.
Due to the considerable size of the literature, we prefer to quote 
the detailed overview by Miranville \cite{mir-CH} and the references therein
for completeness.
Let us only point out the contributions 
\cite{col-gil-spr, cher-mir-zel, gil-mir-sch} dealing with 
well--posedness and
\cite{col-far-hass-gil-spr, col-gil-spr-contr, col-gil-spr-contr2, hinter-weg}
in the direction of distributed and boundary control problems.
Possible relaxations and asymptotics of the 
Cahn--Hilliard equation have been recently studied 
in \cite{bct, bcst1, bcst2, col-scar, scar-VCHDBC}
also with nonlinear viscosity terms.

In the stochastic case, the original contribution dealing with 
Cahn--Hilliard equation is \cite{daprato-deb}, on the 
existence of mild solutions in the case of polynomial potentials.
Further studies have been then carried out in
the works \cite{corn, elez-mike} again in the polynomial setting,
and in \cite{scar-SCH, scar-SVCH} in the case of 
more general potentials in variational framework.
The stochastic Cahn--Hilliard equation with logarithmic potential has been
studied in \cite{deb-zamb, deb-goud, goud}
in relation with reflection measures, and in
\cite{scar-SCHDM} in the case of degenerate mobility.
In the context of phase--field modelling with stochastic 
forcing, it is worthwhile mentioning the contributions
\cite{ant-kar-mill,feir-petc,feir-petc2}, as well as
\cite{bauz-bon-leb, bert, orr-scar} on 
the stochastic Allen-Cahn equation. 
In the direction of optimal control, 
we point out 
\cite{scar-OCSCH} dealing with a distributed optimal 
control problem of the stochastic Cahn--Hilliard equation, 
and the recent work 
\cite{orr-roc-scar} on a stochastic phase--field model 
for tumour growth.

Concerning specifically the Cahn--Hilliard equation with convection,
in the deterministic case well--posedness has been studied in \cite{col-gil-spr-conv}
under general choices of dynamic boundary conditions,
in \cite{della-grass} in a local version with reaction terms, 
while some related optimal velocity control problems have 
been analysed in \cite{col-gil-spr-conv2, col-gil-spr-conv3, roc-spr, zhao-liu, zhao-liu2}.
Also, the relationship between the 
behaviour of the convection term and phase--separation
has been analysed in the recent work \cite{MFFIT-CH}: here, the authors show
that if the velocity field is sufficiently mixing, then no
phase--separation occurs, and the solutions of the 
respective advective Cahn--Hilliard equation converge exponentially
to a homogenous mixed state instead. This may have 
important connections to related optimal control problems
with a target distribution at a final time: in particular, 
the above--mentioned result makes the optimisation problem 
meaningful also when the final target state is not necessarily separated, but is
a homogenous mixed state. Also, it points out how powerful 
the action of the convection term is on the phase--separation,
and motivates the study of phase--optimisation problems where the control
is the velocity itself. 
The convective Cahn--Hilliard equation has also been 
considered in coupled systems, with a further equation equation 
for the velocity field: it is the case, for example, 
of Cahn--Hilliard--Navier--Stokes systems,
studied in \cite{abels-CHNS, frig-grass-spr, frig-grass-spr2, frig-roc-spr}.
By contrast, despite its strong relevance in 
application to stochastic optimal velocity control, 
the convective Cahn--Hilliard has not been analysed yet.
The only results available in the stochastic setting deal with
coupled systems, for example in the context of 
stochastic Cahn--Hilliard--Navier--Stokes models 
\cite{deg-tac-SCHNS,deg-tac-SCHNS2,tac-SCHNS}.
This paper constitutes 
a first contribution to
optimal velocity control for the
stochastic convective Cahn--Hilliard equation.

The literature on stochastic optimal control is also quite extensive:
for a general overview we refer to the monograph \cite{yong-zhou}.
Stochastic optimal control is also studied in 
\cite{fuhr-hu-tess, fuhr-hu-tess2, fuhr-hu-tess3, fuhr-orr, gua-mas-orr}
in the context of the heat equation and reaction-diffusion systems.
For completeness, we refer also to the works
\cite{du-meng, lu-zhang} concerning the stochastic maximal principle.
Relaxation of the optimality conditions have been addressed 
in \cite{brzez-serr} and \cite{barbu-rock-contr}
for dissipative SDPEs and the Schr\"odinger equation, respectively.
Deterministic optimal control problems of stochastic 
reaction--diffusion equations have been analysed in \cite{stann:det_cont}.

Let us describe now the main points that will be addressed in this work.
First of all, we concentrate on the well--posedness of the state--system
\eqref{eq1}--\eqref{eq4}, where the control $\bu$ is arbitrary but fixed.
Using a Yosida approximation on the nonlinearity and 
a time--regularisation on the velocity field, 
we show existence--uniqueness of solutions
by means of variational techniques and stochastic 
compactness arguments. Thanks to monotone analysis 
tools, we are able to cover very general potentials, not 
necessarily of polynomial growth.
Also, we prove
continuous dependence of the variables 
with respect to the control, and
this allows to define a suitable control--to--state map $S:\bu\mapsto(\varphi,\mu)$.
Secondly, we focus on
the optimisation problem, which consists in
minimising a tracking--type cost functional in the form
\[
  J(\varphi,\bu):=
  \frac{\alpha_1}{2}\E\int_Q|\varphi-\varphi_Q|^2
  +\frac{\alpha_2}{2}\E\int_\OO|\varphi(T) - \varphi_T|^2
  +\frac{\alpha_3}{2}\E\int_Q|\bu|^2
\]
subject to the state--system \eqref{eq1}--\eqref{eq4}
and the constraint that $\bu$ is an admissible control, 
meaning that $\bu\in\mU_{ad}$ with $\mU_{ad}$
being a suitable bounded, closed subset of the space
$p$--integrable progressively measurable process
with values in $L^3(\OO)^d$. Here, $\varphi_Q$ and $\varphi_T$
represent some running and final targets, while 
$\alpha_1, \alpha_2, \alpha_3$ are nonnegative weights.

Cost functionals in this form arise very naturally from applications.
Roughly speaking, the optimisation problem amounts to 
identify the optimal way of stirring and mixing the fluid in such a way that
the state variable $\varphi$ is as close as possible to the 
running target $\varphi_Q$ during the evolution and to the final target $\varphi_T$
at the end of the evolution, without wasting too much energy 
in inducing the flow $\bu$. As we have anticipated above, 
a typical example that we have in mind appears in the solidification process of silicon crystals 
in the context of industrial photovoltaic applications \cite{kud, roc-spr}.
Here, a certain mixture of impurities
needs to be moved by convection from within the silicon melt
to its boundary, in order to refine the quality of the final silicon block.
The flow $\bu$ of the fluid behaves then as a control on the silicon melt in order to 
make the relative distribution of impurities $\varphi$ be
close enough to some prescribed targets.
In particular, the final target distribution $\varphi_T$ of impurities 
can be seen here as concentrated on the boundary
and diluted in the interior. 
Analogous applications arise more generally in optimal distribution problems 
of melting materials: the local distribution of some substance contained
in the separating fluid is optimised close to some desired targets 
by inducing a flow in the material itself.

The starting point in the analysis
consists in addressing existence of optimal controls.
This is one of the main differences with respect to the 
deterministic optimal control problem. Indeed, in 
the deterministic setting 
existence of optimal controls follows with no 
particular effort 
from the direct method of calculus of variations, 
since one is able to obtain enough compactness 
from the well--posedness of the state system and the 
boundedness of the set of admissible controls. 
By contrast, in the stochastic 
case these uniform estimates on the minimising sequence of 
controls do not ensure enough compactness in probability,
due to the stochastic nature of the problem itself. Also, 
classical stochastic tools that are usually employed to bypass this problem,
such as the well--known criterion \`a la Gy\"ongy--Krylov, do not work here:
this is due to the non--uniqueness of optimal controls, which is 
caused by the highly nonlinear nature of the minimisation problem.
To overcome this issue, we propose instead a relaxed notion
of optimality, which may be considered as optimality in law, i.e.~requiring 
that the stochastic basis and the Wiener process are part of
the definition of optimal control themselves. This technique mimics 
the definition of probabilistically weak solution for stochastic evolution equations,
and has been employed in other settings such as 
\cite{barbu-rock-contr, orr-roc-scar}.
In this framework, we prove existence of relaxed optimal controls, and we 
show that when one restricts the attention only to deterministic controls
then it is possible to get existence in the classical (probabilistically strong) sense.

We move then to the study of the differentiability properties of 
the control--to--state map $S$. More specifically, we prove that $S$
is G\^ateaux and Fr\'echet differentiable between suitable Banach spaces.
This is done by showing well--posedness of the so--called linearised system,
obtained from \eqref{eq1}--\eqref{eq4} formally differentiating with respect to $\bu$,
and by carefully proving that the unique linearised solution actually coincides with 
the derivative of $S$.
This will allow to explicitly characterise, thanks to the chain rule
in Banach spaces, the derivative of the reduced cost functional $J\circ S$, 
so that the optimisation problem could be seen 
only in terms of the control $\bu$.
Consequently, it is possible to obtain a first rudimental version of necessary conditions
for optimality, by imposing the classical first--order 
variational inequality $D(J\circ S)(\bu)\geq0$ on a given optimal control.

The last part of the paper aims at refining the first version
of necessary conditions, by removing any explicit dependence 
on the linearised variables. This is done by introducing and studying 
a suitable adjoint problem, which is formally related to the dual 
problem of the linearised system.
The adjoint problem consists of a backward--in--time stochastic 
partial differential equation, and its analysis is the 
most challenging point of the work.
The first main difficulty is indeed the backward nature of the equation:
although this is not a great limitation in deterministic problems, 
in the stochastic case it calls for the introduction of 
an extra variable, in order to preserve adaptability of the processes in play,
and requires different analytical techniques such as 
martingale representation theorems.
The second and most crucial difficulty depends instead on the 
nonlinear nature of the system.
Indeed, the presence of the nonlinear term $\Psi''(\varphi)$
and the dual structure of the equation
prevent from obtaining 
uniform estimates directly on the adjoint system. 
Consequently, well--posedness cannot be
obtained classically by
tackling the adjoint problem straightaway, and a different idea is needed.
In this regard, we use a duality method.
We consider a more general version of the linerised system, 
where an arbitrary forcing term is added, and we show that 
this is well--posed and the solutions depend continuously 
on the forcing term. Then, we prove that such system is
in duality with the adjoint problem that we want to study, and
this allows to recover by comparison some first uniform 
estimates on the adjoint variables. 
This tool is extremely powerful, as it allows to 
bound the adjoint variables without even
working on the adjoint system itself: the main intuition 
behind this is that the linearised system is usually much simpler to study,
and the duality between linerised--adjoint systems 
allows to ``transfer'' uniform bounds on the solutions
from one problem to the other.
Once these first crucial
estimates are obtained, using 
classical techniques we are then able to prove well--posedness 
of the adjoint problem. Lastly, the duality relation 
is employed to refine the first--order conditions for optimality
and to write them as a variational inequality only depending on
the intrinsic adjoint variables.

The main novelty of the work is the presence of two sources of 
randomness in equation \eqref{eq1}, accounting for 
noises both in the phase--separation process and in the 
flow--inducing procedure. 
As interesting as it may be 
from the applied point of view, 
certainly this novel framework does not
come without effort on the mathematical side.
Indeed, let us stress that the fact that $\bu$ is assumed to be a 
stochastic process, and not a deterministic function, causes
several non--trivial issues in estimating the solutions:
this is due to a lack of satisfactory computational tools 
of Gronwall--type in the genuinely pure stochastic case.
Such difficulties are evident especially in the study of the forward
problems, i.e.~in the state system \eqref{eq1}--\eqref{eq4} and in the 
corresponding linearised system. Here, the idea 
is to argue instead combining carefully the H\"older inequality 
and several iterative patching arguments, in order to avoid 
applying the Gronwall lemma, which does not work.
In the adjoint problem, the situation is slightly better:
we will show that the backward nature of the equation allows indeed to 
use a very general and recent backward--in--time version of the stochastic Gronwall lemma
(see Lemma~\ref{lem:gronwall} below).

We conclude by summarising here the structure of the paper.
Section~\ref{sec:main} contains the description of the setting
of the work, the precise assumptions, and the main results that we prove.
In Section~\ref{sec:WP} we prove well--posedness of the state--system,
while Section~\ref{sec:ex_optt} focuses on the existence of optimal controls.
Then, in Sections~\ref{sec:lin} and \ref{sec:ad} we study the 
linearised system and the adjoint system, respectively.
Finally, in Section~\ref{sec:nec} we prove the two versions of 
first--order conditions for optimality.

%%%%%%%%%%%%%%%%%%%%%%%%%%%%%%%%%%%%%%%%%%%%%%%%

\section{Setting and assumptions}
\label{sec:main}
In this section we specify the general setting, notation, 
and assumptions of the work. We then 
present the main results of the paper.

Let $(\Omega,\cF, (\cF_t)_{t\in[0,T]}, \P)$ be a filtered probability space
satisfying the usual conditions, where $T>0$ is a fixed final time and
$W$ is a cylindrical Wiener process on a separable Hilbert space $K$.
For convenience, let us fix now once and for all 
a complete orthonormal system $(e_j)_j$ of $K$.
The progressive $\sigma$-algebra on $\Omega\times[0,T]$
is denoted by $\cP$.

As far as notation is concerned, the dual
of a given real Banach space $E$ is denoted by $E^*$,
and the duality pairing between $E^*$ and $E$ is denoted by 
$\ip{\cdot}{\cdot}_{E^*,E}$.
Weak convergence in $E$ and weak$^*$ convergence in $E^*$
will be denoted by the respective symbols $\wto$ and $\wstarto$.
Also, for all $q\in[1,+\infty]$
we employ the usual symbols $L^q(\Omega; E)$
and $L^q(0,T; E)$ for the spaces
of $q$--Bochner integrable functions, and 
$C^0([0,T]; E)$ and $C^0_w([0,T]; E)$
for the spaces of strongly and weakly continuous 
functions from $[0,T]$ to $E$, respectively.
For spaces of stochastic processes,
we use the notation $L^{q_1}_\cP(\Omega; L^{q_2}(0,T; E))$
to further specify that measurability is also intended 
with respect to the progressive $\sigma$-algebra $\cP$.
In the case that $q>1$ and $E$ is separable, 
we explicitly set 
$L^q_w(\Omega; L^\infty(0,T; E^*))$ as the dual space
of $L^{\frac{q}{q-1}}(\Omega; L^1(0,T; E))$,
which we recall can be characterised \cite[Thm.~8.20.3]{edwards} as 
the space of weak*-measurable 
random variables $y:\Omega\to L^\infty(0,T; E^*)$
with finite $q$-moment in $\Omega$.
Finally, if $E_1$ and $E_2$ are separable Hilbert spaces, 
we use the notation $\cL^2(E_1,E_2)$ for the space of 
Hilbert-Schmidt operators from $E_1$ to $E_2$.

In the proofs, the symbol $c$ is reserved to denote 
any generic positive constant, whose value depends 
on the structure of the problem and may be updated 
from line to line in the proofs.

Let $\OO\subset\erre^d$ ($d\geq2$) be a smooth bounded domain.
We use the classical notation $Q:=(0,T)\times\OO$, 
$Q_t:=(0,t)\times\OO$, and $Q_t^T:=(t,T)\times\OO$
for every $t\in(0,T)$. The outward normal unit vector on the boundary $\partial\OO$
is denoted by ${\bf n}$.
We introduce the functional spaces
\begin{align*}
  &H:=L^2(\OO)\,,  \qquad
  V_1:=H^1(\OO)\,, \\
  &V_2:=\{v\in H^2(\OO):\;{\bf n}\cdot\nabla v = 0\quad\text{a.e.~on } \partial\OO\}\,,
  \qquad V_3:=V_2\cap H^3(\OO)\,,
\end{align*}
endowed with their natural norms $\norm{\cdot}_H$, 
$\norm{\cdot}_{V_1}$, $\norm{\cdot}_{V_2}$, and $\norm{\cdot}_{V_3}$
respectively. We identify $H$ to its dual, so that we have the
continuous and dense inclusions
\[
  V_3\embed V_2\embed V_1\embed H \embed V_1^*\,.
\]
For all $y\in V_1^*$ we use the notation $y_\OO:=\frac1{|\OO|}\ip{y}{1}$
for the spatial mean of $y$, and define the subspaces of zero-mean elements as
\[
  V_{1,0}^*:=\{y\in V_1^*:\;y_\OO=0\}\,, \qquad  H_0:=H\cap V_{1,0}^*\,, \qquad
  V_{1,0}:=V_1\cap H_0\,.
\]
Let us  recall that the variational formulation of the Laplace operator 
with Neumann conditions
\[
  \mathcal L: V_1\to V_1^*\,, \qquad
  \ip{\mathcal L y}{\zeta}:=\int_\OO\nabla y\cdot\nabla\zeta\,,\quad y,\zeta\in V_1\,,
\]
is a well-defined linear operator, and its restriction to $V_{1,0}$
is an isomorphism onto the space $V_{1,0}^*$. 
Its inverse $\mN: V_{1,0}^*\to V_{1,0}$
is the resolvent operator associated to the abstract elliptic problem on $\OO$
with homogenous Neumann conditions, meaning that for all 
$y\in V_{1,0}^*$ the element $z:=\mN y\in V_{1,0}$ is the unique solution
with null mean to
\[
  \begin{cases}
  -\Delta z = y\quad&\text{in } \OO\,,\\
  \partial_{\bf n}z=0\quad&\text{in } \partial\OO\,.
  \end{cases}
\]
As a consequence of the Poincar\'e--Wirtinger inequality,
it is immediate to check that 
\[
  \zeta\mapsto \norm{\nabla\mN(\zeta-\zeta_\OO)}_H^2 + |\zeta_\OO|^2\,,\qquad
  \zeta\in V_1^*\,,
\]
yields an equivalent norm on $V_1^*$. In particular, it follows the compactness inequality 
\beq
  \label{comp_ineq}
  \forall\,\eps>0\,,\quad\exists\,c_\eps>0:\quad
  \norm{y}_H^2\leq\eps\norm{\nabla y}_{H}^2 + c_\eps\norm{\nabla\mathcal Ny}_H^2
  \quad\forall\,y\in V_{1,0}\,.
\eeq

We introduce the space
  \[
  U:=\left\{{\bf u}\in L^3(\OO):\quad
  \div\bu=0\,,\quad\bu\cdot{\bf n}=0 \text{ a.e.~on } \partial\OO\right\}\,,
  \]
  where the divergence is intended in the sense of distributions on $\OO$.
  The space of velocity controls $\bf u$ that we focus on will be 
  \[
  \mathcal U:=L^\infty_\cP(\Omega; L^p(0,T; U))\,, \qquad p\in(2,+\infty)\,.
  \]
  Let us note that this includes as a special case the choice of deterministic controls,
  which has also received a strong mathematical interest on its own: 
  see for instance Stannat \& Wessels \cite{stann:det_cont}. Indeed, 
  we can set
  \[
  \mathcal U^{det}:=L^p(0,T; U)\subset \mU\,.
  \]

The following assumptions on the problem
will be in force throughout the paper.
\begin{description}
  \item[A1] $\Psi:\erre\to\erre$ is of class $C^2$, $\Psi'(0)=0$, and there exist
  $C_\Psi>0$ and $\gamma\in[1,2]$ such that 
  \begin{align*}
  \Psi''(r)\geq -C_\Psi \qquad&\forall\,r\in\erre\,,\\
  |\Psi'(r)| +  |\Psi''(r)|^\gamma\leq C_\Psi(1+ \Psi(r)) \qquad&\forall\,r\in\erre\,.
  \end{align*}
  Let us point out that the classical polynomial double--well potential $\Psi_{pol}$
  satisfies these assumptions with $\gamma=2$.
  Nonetheless, by allowing also the smaller values $\gamma\in[1,2]$
  we are able to include possibly more singular potential, such as
  first--order exponentials.
  We set $\beta:r\mapsto \Psi'(r) + C_\Psi r$, $r\in\erre$:
  then $\beta:\erre\to\erre$ is a $C^2$ nondecreasing function, hence it can
  be identified with a maximal monotone (single-valued) graph in $\erre\times\erre$.
  Let us also denote by $\widehat\beta:\erre\to[0,+\infty)$ the convex lower semicontinuous
  function with $\widehat\beta(0)=0$. 
  \item[A2] $\varphi_0 \in V_1$ and $\Psi(\varphi_0)\in L^1(\OO)$.
  \item[A3] $B:V_1\to\cL^2(K,V_1)$ and there exists a constant $C_B>0$ such that 
  \begin{align*}
    \norm{B(y_1)-B(y_2)}_{\cL^2(K,H)} \leq C_B\norm{y_1-y_2}_{H} 
    &\quad\forall\,y_1,y_2\in H\,,\\
    \norm{B(y)}_{\cL^2(K,V_1)} \leq C_B\left(1+\norm{y}_{V_1} \right)
    &\quad\forall\,y\in V_1\,,\\
    \sum_{j=0}^\infty
    \norm{B(y)e_j}_{L^{\frac{2\gamma}{\gamma-1}}(\OO)}^2
    \leq C_B
    &\quad\forall\,y\in H\,.
  \end{align*}
  Moreover, we prescribe that
  \[
    B:V_1\to\cL^2(K,V_{1,0}) \qquad\text{in case of multiplicative noise}\,.
  \]
  Let us note that in case of additive noise $B\in\cL^2(K,V_1)$, these conditions 
  are trivially satisfied for all $\gamma\in(1,2]$ if $d=2$ 
  and for all $\gamma\in[3/2,2]$ if $d=3$: in particular, the classical 
  polynomial case in dimension two and three is always covered.
  In the genuine multiplicative noise case,
  i.e.~when $B$ is not constant in $V_1$, we also suppose that 
  $B$ is $\cL^2(K,V_{1,0})$-valued: 
  this amounts to requiring that 
  the noise is conservative, in the sense that it preserves the mean
  $\varphi_\OO$ of the phase-variable. A direct consequence 
  is the conservation of mass, which is a fundamental feature of 
  Cahn--Hilliard-type evolutions. 
  This hypothesis on the noise is very classical and natural in literature:
  for example, let us stress that a relevant multiplicative choice of $B$ can be given as
  \[
  B(y)e_j:=h_j(y) - (h_j(y))_\OO\,, \quad y\in V_1\,,\quad j\in\enne\,,
  \]
  where the sequence $(h_j)_j\subset W^{1,\infty}(\erre)$ is such that
  \[
  C_B^2:=\sum_{j=0}^\infty\norm{h_j}^2_{W^{1,\infty}(\erre)}<+\infty\,.
  \]
  It is not difficult to show that this example allows for all values of $\gamma\in[1,2]$
  in every space-dimension $d=2,3$.
\end{description}

In the context of the optimal velocity control, it will be useful to 
introduce a polynomial-growth assumption on $\Psi$.
This will be necessary only in the study of the optimisation problem, but 
is not needed for the well--posedness of the state system.
\begin{description}
  \item[C1] it holds that $\gamma=2$ in {\bf A1} and
  \[
  |\Psi''(r)|\leq C_\Psi(1+|r|^2) \qquad\forall\,r\in\erre\,.
  \]
  Such requirement is very natural in the Cahn--Hilliard context,
  since it satisfied by the classical choice of 
  the polynomial double-well potential $\Psi_{pol}$ of degree $4$.
\end{description}

The first main result of the paper states existence and uniqueness
of strong solutions, and their continuous dependence with respect to the velocity field.
\begin{thm}
  \label{thm1}
  Assume {\bf A1--A3}. Then, for every $\bu\in\mathcal U$, there exists a unique 
  pair $(\varphi,\mu)$ with 
  \begin{align*}
  &\varphi \in L^p_\cP\left(\Omega; W^{s,p}(0,T; V_1^*)\cap  C^0([0,T]; H)
  \cap L^2(0,T; V_2)\right) \cap L^p_w(\Omega;L^\infty(0,T; V_1))\,,\\
  &\mu = -\Delta\varphi + \Psi'(\varphi)\in L^{p/2}_\cP(\Omega; L^2(0,T; V_1))\,,
  \end{align*}
  for all $s\in(0,1/2)$, and such that
  \begin{align*}
  (\varphi(t),\zeta)_H + \int_{Q_t}\nabla\varphi\cdot\nabla\zeta 
  -\int_{Q_t} \varphi\bu\cdot\nabla\zeta= 
  (\varphi_0,\zeta)_H + \left(\int_0^tB(\varphi(s))\,\d W(s), \zeta\right)_H
  \qquad\forall\,\zeta\in V_1\,,
  \end{align*}
  for every $t\in[0,T]$, $\P$-almost surely. Furthermore, 
  there exists 
  a constant $K>0$, only depending on the structure of the problem,
  such that for all $\bu\in\mathcal U$, the respective solution $(\varphi,\mu)$
  satisfies 
  \begin{align}
    \nonumber
    &\norm{\varphi}_{L^p(\Omega; L^\infty(0,T; V_1))\cap L^p_\cP(\Omega; L^2(0,T; V_2))}
    +\norm{\mu}_{L^{p/2}_\cP(\Omega; L^2(0,T; V_1))}
    +\norm{\Psi(\varphi)}_{L^{p/2}(\Omega;L^\infty(0,T; L^1(\OO)))}\\
    &\qquad+\norm{\Psi'(\varphi)}_{L^{p/2}_\cP(\Omega; L^2(0,T; H))} +
    \norm{\Psi''(\varphi)}_{L^{\gamma p/2}(\Omega;L^\infty(0,T; L^\gamma(\OO)))}
    \label{bound}
    \leq K\left[1 + \norm{\bu}_{\mU}^{\frac{2}{p-2}}\right]\,,
  \end{align}
  and for every $\{\bu_i\}_{i=1,2}\subset\mathcal U$,
  the respective solutions $\{(\varphi_i, \mu_i)\}_{i=1,2}$ verify
  \begin{align}
  &\norm{\varphi_1-\varphi_2}_{L^p_\cP(\Omega; C^0([0,T]; V_1^*)\cap 
   L^2(0,T; V_1))}
   \label{cont}
   \leq K \left[1+\norm{\bu_1}_{\mU}^{\frac{2}{p-2}}\right]
   \left[1+\norm{\bu_2}_{\mU}^{\frac{2}{p-2}}\right]
  \norm{\bu_1-\bu_2}_{\mU}\,.
  \end{align}
  Lastly, if also {\bf C1} holds, then 
  \begin{align}
  \nonumber
  &\norm{\varphi_1-\varphi_2}_{L^{p/3}_\cP(\Omega; C^0([0,T]; H)\cap 
   L^2(0,T; V_2))} + \norm{\mu_1-\mu_2}_{L^{p/3}_\cP(\Omega; L^2(0,T; H))}\\
   \label{cont2}
   &\leq K\left[1 + \norm{\bu_1}_{\mU}^{\frac{4}{p-2}}
   + \norm{\bu_2}_{\mU}^{\frac{4}{p-2}}\right]
   \left[1+\norm{\bu_1}_{\mU}^{\frac{2}{p-2}}\right]
   \left[1+\norm{\bu_2}_{\mU}^{\frac{2}{p-2}}\right]
  \norm{\bu_1-\bu_2}_{\mU}\,.
  \end{align}
\end{thm}

Once the analysis of well-posedness of the state-system has been addressed, 
we can turn our attention to the optimal velocity control problem.
As far as the controls are concerned, we consider classical box-contraints 
on the velocity controls, by defining the set of admissible controls as
\[
  \mathcal U_{ad}:=\left\{\bu \in\mathcal U:\;\norm{\bu}_{L^p(0,T; U)} \leq L
  \quad\P\text{-a.s.}\right\}\,,
\]
where $L>0$ is a prescribed constant. The prescription
of a box--constraint on the admissible controls is classical on the mathematical 
side. In applications, the constant $L$ is typically related to the 
maximum capacity of the flow--inducing devices that convey the velocity field.
It will be useful to introduce 
an enlarged bounded open set $\widetilde{\mathcal U}_{ad}$ in 
$\mU$ containing $\mathcal U_{ad}$, as
\[
  \widetilde{\mathcal U}_{ad}:=\left\{\bu \in\mathcal U:\;\norm{\bu}_{\mU} < L+1\right\}\,.
\]
Analogously, we introduce the corresponding spaces of admissible deterministic controls as
\[
  \mU_{ad}^{det}:=\mU^{det}\cap\mathcal U_{ad}\,,\qquad
  \widetilde\mU_{ad}^{det}:=\mU^{det}\cap \widetilde{\mathcal U}_{ad}\,.
\]
The cost functional that we study is 
of quadratic tracking-type and reads
\begin{align}
  \nonumber
  &J:L^2_\cP(\Omega;C^0([0,T]; H))\times L^2_\cP(\Omega; L^2(0,T; H^d)) \to \erre\,,\\
  \label{cost}
  &J(\varphi,\bu):=\frac{\alpha_1}{2}\E\int_Q|\varphi-\varphi_Q|^2
  +\frac{\alpha_2}{2}\E\int_\OO|\varphi(T) - \varphi_T|^2
  +\frac{\alpha_3}{2}\E\int_Q|\bu|^2\,,\\ 
  \nonumber
  &\qquad\qquad\;\;(\varphi,\bu)\in L^2_\cP(\Omega;C^0([0,T]; H))\times 
  L^2_\cP(\Omega; L^2(0,T; H^d)) \,,
\end{align}
where $\alpha_1,\alpha_2,\alpha_3$ are nonnegative constants with
$\alpha_1+\alpha_2+\alpha_3>0$ and the targets are fixed with
\[
  \varphi_Q \in L^2_\cP(\Omega; L^2(0,T; H))\,, \qquad 
  \alpha_2\varphi_T\in L^2(\Omega,\cF_T; H)\,.
\]
The optimal velocity control consists in the following:
\begin{description}
  \item[(CP)] minimise the cost functional $J$
  with the constraints that $\bu$ belongs to $\mathcal U_{ad}$ and 
  $\varphi$ is the unique corresponding solution component to 
  the state system \eqref{eq1}--\eqref{eq4}.
\end{description}
By virtue of the well-posedness Theorem~\ref{thm1} it is well-defined
the {\em control-to-state map}
\[
  S:\widetilde\mU_{ad}\to 
  \left[L^p_\cP\left(\Omega; C^0([0,T]; H)
  \cap L^2(0,T; V_2)\right) \cap L^p_w(\Omega;L^\infty(0,T; V_1))\right]\times
  L^{p/2}_\cP(\Omega; L^2(0,T; V_1))
\]
as
\[
  S(\bu)=(S_1(\bu), S_2(\bu)):=(\varphi,\mu)\,, \qquad\bu\in\widetilde\mU_{ad}\,.
\]
This implies that the optimal control problem can be reduced 
to the only variable ${\bu}$, by introducing 
the so-called {\em reduced} cost functional as
\[
  \widetilde J:\widetilde \mU_{ad}\to \erre\,, \qquad
  \widetilde J(\bu):=J(S_1(\bu), \bu)\,, \quad \bu\in\widetilde\mU_{ad}\,.
\]
\begin{remark}
Clearly the well--posedness result in
Theorem~\ref{thm1} continues to hold on any new stochastic basis
$(\Omega', \cF', \P', W')$,
provided to analogously define the new spaces of 
controls $\mU'$, $\mU'_{ad}$, and $\widetilde\mU_{ad}'$.
Hence, if also $(\varphi_Q',\varphi_T')$ are some 
new targets on $(\Omega', \cF', \P')$ with the same law
of $(\varphi_Q, \varphi_T)$, one can define the corresponding 
cost functional $J'$, the corresponding control--to--state map $S'$,
and the new reduced cost functional $\widetilde J'$ on the new probability space,
by simply replacing $\Omega$ with $\Omega'$.
\end{remark}

With this notations, we can state the exact definition of optimal control as follows.
As anticipated, we also give some relaxed notions of optimality,
one based on the concept of optimality--in--law and the other 
obtained minimising only on the deterministic controls.
\begin{defin}
  \label{def:OC}
  An optimal control for {\bf(CP)} is an element $\bu\in\mU_{ad}$ such that
  \[
  \widetilde J(\bu)= \inf_{\bv\in\mU_{ad}}\widetilde J(\bv)\,.
  \]
  A relaxed optimal control for {\bf(CP)} is a family 
  $\left(\Omega', \cF', (\cF'_t)_{t\in[0,T]}, \P', W',\varphi_Q', \varphi_T',\bu'\right)$
  where $(\Omega', \cF', \P')$ is a probability space,
  $(\cF'_t)_{t\in[0,T]}$ is a filtration satisfying the usual conditions,
  $W'$ is a $K$-cylindrical Wiener process on it, 
  $\alpha_1\varphi_Q'\in L^2_\cP(\Omega'; L^2(0,T; H))$ and 
  $\alpha_2\varphi_T'\in L^2(\Omega',\cF_T'; H)$
  have the same laws of $\alpha_1\varphi_Q$ and $\alpha_2\varphi_T$, 
  respectively, and $\bu'\in\mU_{ad}'$ satisfies 
  \[
  \widetilde J'(\bu')= \inf_{\bv\in\mU_{ad}}\widetilde J(\bv)\,.
  \]
  A deterministic optimal control for {\bf(CP)} is an element $\bu\in\mU_{ad}^{det}$ such that
  \[
  \widetilde J(\bu)= \inf_{\bv\in\mU_{ad}^{det}}\widetilde J(\bv)\,.
  \]
\end{defin}

Our first result in the analysis of the optimisation problem {\bf(CP)}
concerns existence optimal controls. 
It is worthwhile noting that due to the non-uniqueness of optimal controls,
in the genuinely stochastic case one can only
show existence of relaxed optimal controls:
this is typical in highly nonlinear stochastic optimal control problems,
se for example \cite{barbu-rock-contr, scar-OCSCH}.
By contrast, we 
show that deterministic optimal controls always exist.
\begin{thm}
  \label{thm2}
  Assume {\bf A1--A3}. Then, there exist a relaxed optimal control $\bu$
  and a deterministic optimal control $\bu^{det}$ for problem {\bf (CP)}.
\end{thm}

Once existence of minimisers for {\bf(CP)} is proved, 
we can now turn to the main focus of the work, i.e~the investigation 
of necessary conditions for optimality.
The first main step in this direction is the study 
of the differentiability of the control-to-state map $S$,
along with the characterisation of its derivative through 
the analysis of the linearised state system. 
This will allow to obtain a first version of first-order conditions for optimality 
by means of a suitable variational inequality involving the derivative 
of the reduced cost functional.
In this direction, we introduce the assumptions
\begin{description}
\item[C2] the map $B:V_1\to\cL^2(K,H)$ is of class $C^1$.
  Let us point out that this implies together with {\bf A3}
  that $\norm{DB(y)\zeta}_{\cL^2(K,H)}\leq C_B\norm{\zeta}_H$
  for all $y,\zeta\in V_1$.
  Moreover, let us stress this requirement is very natural, 
  and it is satisfied for instance in the relevant
  example described in {\bf A3}, provided to 
  replace $W^{1,\infty}(\erre)$ with $W^{1,\infty}(\erre)\cap C^1(\erre)$.
\item[C3] $\Psi$ is of class $C^3$, $DB\in C^{0,1}(V_1; \cL(V_1,\cL^2(K,H)))$, and it holds that 
  \[
  |\Psi'''(r)|\leq C_\Psi(1+|r|)\qquad\forall\,r\in\erre\,.
  \]
  This is a refinement of assumptions {\bf C1--C2} and ensures, as we will see, 
  better differentiability properties for $S$. Still, {\bf C3} is satisfied by the polynomial 
  potential $\Psi_{pol}$ and the relevant noise coefficient 
  described in {\bf A3}, provided to
  replace $W^{1,\infty}(\erre)$ with $W^{2,\infty}(\erre)$.
\end{description}
The linearised system
can be formally obtained by differentiating the state system \eqref{eq1}--\eqref{eq4}
with respect to the control $\bu$ in a given direction $\bh\in\mU$, and reads
\begin{align}
    \label{eq1_lin}
    \d\theta_\bh - \Delta \nu_\bh\,\d t + \bh\cdot\nabla\varphi\,\d t + \bu\cdot\nabla\theta_\bh\,\d t
    = DB(\varphi)\theta_\bh\,\d W \qquad&\text{in } (0,T)\times \OO\,,\\
    \label{eq2_lin}
    \nu_\bh=-\Delta \theta_\bh + \Psi''(\varphi)\theta_\bh \qquad&\text{in } (0,T)\times \OO\,,\\
    \label{eq3_lin}
    {\bf n}\cdot\nabla\theta_\bh = {\bf n}\cdot\nabla\nu_\bh = 0 
    \qquad&\text{in } (0,T)\times\partial \OO\,,\\
    \label{eq4_lin}
     \theta_\bh(0)=0 \qquad&\text{in } \OO\,.
\end{align}
The next result ensures exactly that the linearised system \eqref{eq1_lin}--\eqref{eq4_lin}
is well-posed in a suitable variational sense, and that the unique solution
to \eqref{eq1_lin}--\eqref{eq4_lin} coincides with the 
derivative of the control-to-state map $S$ in the point $\bu$ along the direction $\bh$.
\begin{thm}
  \label{thm3}
  Assume {\bf A1--A3}, {\bf C1--C2}, and $p>3$. 
  Then, for all $\bu\in \widetilde\mU_{ad}$ and
  $\bh\in\mU$, setting $\varphi:=S_1(\bu)$,
  there exists a unique pair $(\theta_\bh,\nu_\bh)$ with
  \begin{align*}
  &\theta_\bh \in 
  L^{p}_{\cP}\left(\Omega; C^0([0,T]; V_1^*)\cap L^2(0,T; V_1)\right)
  \cap L^{p/3}_\cP\left(\Omega; C^0([0,T]; H)\cap  L^2(0,T; V_2)\right)\,,\\
  &\nu_\bh= -\Delta \theta_\bh + \Psi''(\varphi)\theta_\bh\in L^{p/3}_{\cP}(\Omega; L^2(0,T; H))\,,
  \end{align*}
  such that, for every $t\in[0,T]$, $\P$-almost surely,
  \begin{align*}
    (\theta_\bh(t),\zeta)_H - \int_{Q_t}\nu_\bh\Delta\zeta
  -\int_{Q_t}(\varphi\bh + \theta_\bh\bu)\cdot\nabla\zeta= 
  \left(\int_0^tDB(\varphi(s))\theta_\bh(s)\,\d W(s), \zeta\right)_H
  \quad\forall\,\zeta\in V_2\,.
  \end{align*}
  Furthermore,  the control-to-state map $S_1$
  is G\^ateaux-differentiable in the following sense: 
  for all $\bu\in \widetilde\mU_{ad}$ and $\bh\in\mU$, as $\delta\searrow0$ it holds that
  \begin{align*}
  \frac{S_1(\bu+\delta \bh) - S_1(\bu)}{\delta}\to \theta_\bh \qquad&\text{in } 
  L^\ell_\cP(\Omega; L^2(0,T; V_1))\quad\forall\,\ell\in[1,p)\,,\\
  \frac{S_1(\bu+\delta \bh) - S_1(\bu)}{\delta}\wstarto \theta_\bh \qquad&\text{in } 
  L^{p}_w\left(\Omega; L^\infty(0,T; V_1^*)\right)\cap
  L^{p}_\cP\left(\Omega; L^2(0,T; V_1)\right)\,,\\
  \frac{S_1(\bu+\delta \bh) - S_1(\bu)}{\delta}\wstarto \theta_\bh \qquad&\text{in } 
  L^{p/3}_w\left(\Omega; L^\infty(0,T; H)\right)\cap
  L^{p/3}_\cP\left(\Omega; L^2(0,T; V_2)\right)\,,\\
  \frac{S_1(\bu+\delta \bh)(t) - S_1(\bu)(t)}{\delta}\wto \theta_\bh(t) \qquad&\text{in } 
  L^{p/3}(\Omega, \cF_t; H) \quad\forall\,t\in[0,T]\,.
  \end{align*}
  Moreover, if $p\geq7$ and {\bf C3} holds, then $S_1$ is also Fr\'echet-differentiable as a map
  \[
  S_1:\widetilde\mU_{ad}\to L^{p/7}_\cP(\Omega; C^0([0,T]; V_1^*)\cap L^2(0,T; V_1))\,.
  \]
\end{thm}

The second step in the analysis of necessary conditions for optimality 
consists in studying the so-called adjoint system
and by proving a suitable duality relation with respect 
to the linearised system. The adjoint system can be formally obtained 
as the dual system of \eqref{eq1_lin}--\eqref{eq4_lin}, and reads
\begin{align}
  \nonumber
  -\d P -\Delta \tilde P \,\d t + \Psi''(\varphi)\tilde P\,\d t - \bu\cdot\nabla P\,\d t
  \qquad\qquad& \\
  \label{eq1_ad}
  =\alpha_1(\varphi-\varphi_Q)\,\d t+
  DB(\varphi)^*Z\,\d t  - Z\,\d W
  \qquad&\text{in } (0,T)\times\OO\,,\\
  \label{eq2_ad}
  \tilde P=-\Delta P
  \qquad&\text{in } (0,T)\times\OO\,,\\
  \label{eq3_ad}
  {\bf n}\cdot\nabla P = {\bf n}\cdot\nabla \tilde P = 0
  \qquad&\text{in } (0,T)\times\partial\OO\,,\\
  \label{eq4_ad}
  P(T)=\alpha_2(\varphi(T)-\varphi_T)
  \qquad&\text{in } \OO\,.
\end{align}
Let us point out that the adjoint system is backward in time:
due to the stochastic framework of the problem, 
this necessarily requires the introduction of the additional variable $Z$
in view of the classical martingale representation theorems.
The situation here is then much more complex than the deterministic one:
the variable of the adjoint system is indeed the couple $(P,Z)$,
with $\tilde P$ being an auxiliary variable.
Due to the difficulty of analysis of the adjoint system,
we will need to require more regularity on the targets, namely
\begin{description}
\item[C4] $p\geq6$ and it holds that 
  \[
  \alpha_1\varphi_Q \in L^{\frac{2p}{p-4}}_\cP(\Omega; L^2(0,T; H))\,, \qquad 
  \alpha_2\varphi_T\in L^{\frac{2p}{p-4}}(\Omega,\cF_T; V_1)\,.
  \]
\end{description}
The next result ensures that the adjoint system 
\eqref{eq1_ad}--\eqref{eq4_ad} is well-posed in a suitable variational sense,
and state a duality relation between \eqref{eq1_lin}--\eqref{eq4_lin}
and \eqref{eq1_ad}--\eqref{eq4_ad}.
\begin{thm}
  \label{thm4}
  Assume {\bf A1--A3}, {\bf C1--C2}, and {\bf C4}.
  Then, for all $\bu\in \widetilde\mU_{ad}$,
  setting $\varphi:=S_1(\bu)$,
  there exists a triplet $(P, \tilde P,Z)$, with 
  \begin{align*}
  &P\in L^2_\cP(\Omega; C^0([0,T]; V_1)\cap L^2(0,T; V_3))\,,\\
  &\tilde P=\mathcal L P \in L^2_\cP(\Omega; C^0([0,T]; V_1^*)\cap L^2(0,T; V_1))\,,\\
  &Z \in L^2_\cP(\Omega; L^2(0,T; \cL^2(K,V_1)))\,,
  \end{align*}
  such that, for every $t\in[0,T]$, $\P$-almost surely,
  \begin{align*}
    &\left(P(t), \zeta\right)_H
    +\int_{Q_t^T}\nabla\tilde P\cdot\nabla\zeta
    +\int_{Q_t^T}\Psi''(\varphi)\tilde P\zeta
    +\int_{Q_t^T}P\bu\cdot\nabla\zeta\\
    &\qquad
    =\left(\alpha_2(\varphi(T)-\varphi_T), \zeta\right)_H
    +\int_{Q_t^T}DB(\varphi)^*Z\zeta
    -\left(\int_t^TZ(s)\,\d W(s), \zeta\right)_H
    \qquad\forall\,\zeta\in V_1\,.
  \end{align*}
  Furthermore, the solution components $\nabla P$, $\tilde P$, and $\nabla Z$
  are unique in the spaces $L^2_\cP(\Omega; C^0([0,T]; H^d))$, 
  $L^2_\cP(\Omega; C^0([0,T]; V_1^*))$, and 
  $L^2_\cP(\Omega; L^2(0,T; \cL^2(K,H^d)))$, respectively.
\end{thm}

At this point, we are finally ready to state the necessary conditions 
for optimality: more specifically, we present here two different versions.
The first one is deduced directly by the characterisation of the derivative of $S_1$
in Theorem~\ref{thm3}, and consists of a variational inequality depending also
on the linearised variables. 
The second one is a refinement of this, as it employs the adjoint problem
and only depends on the intrinsic adjoint variables $(P,\tilde P, Z)$,
not on the linearised ones.
\begin{thm}
  \label{thm5}
  Assume {\bf A1--A3}, {\bf C1--C2}, and $p\geq6$.
  If $\bu\in\mathcal U_{ad}$ is an optimal control for 
  {\bf(CP)} and $\varphi:=S_1(\bu)$
  is its respective optimal state,
  then
  \beq\label{NC}
  \alpha_1\E\int_Q(\varphi - \varphi_Q)
  \theta_{\bv-\bu} + \alpha_2\E\int_\OO(\varphi(T)-\varphi_T)\theta_{\bv-\bu}(T)
  +\alpha_3\E\int_Q\bu\cdot(\bv-\bu) \geq 0 \qquad\forall\,\bv\in\mathcal U_{ad}\,,
  \eeq
  where $\theta_{\bv-\bu}$ is the unique first solution component of
  the linearised system \eqref{eq1_lin}--\eqref{eq4_lin}
  with the choice $\bh:=\bv-\bu$, in the 
  sense of Theorem~\ref{thm3}.
\end{thm}

\begin{thm}
\label{thm6}
  Assume {\bf A1--A3}, {\bf C1--C2}, and {\bf C4}.
  If $\bu\in\mathcal U_{ad}$ is an optimal control for 
  {\bf(CP)} and $\varphi:=S_1(\bu)$
  is its respective optimal state,
  then
  \beq
  \label{NC2}
  \E\int_Q(\varphi\nabla P + \alpha_3\bu)\cdot(\bv-\bu) \geq 0 \qquad\forall\,\bv\in\mathcal U_{ad}\,,
  \eeq
  where $\nabla P$ is the uniquely-determined solution component of
  the adjoint system \eqref{eq1_ad}--\eqref{eq4_ad} 
  in the sense of Theorem~\ref{thm4}.
  In particular, if $\alpha_3>0$,
  then $\bu$ is the orthogonal projection of $-\frac1{\alpha_3}\varphi\nabla P$
  on the closed convex set $\mathcal U_{ad}$
  in the Hilbert space $L^2_\cP(\Omega;L^2(0,T; H^d))$.
\end{thm}

\begin{remark}
  Let us comment on the necessary condition for optimality. 
  When handling the optimisation problem in practice, the main role
  of condition \eqref{NC2} is to restrict the class of possible candidates 
  to be optimal controls. Roughly speaking, the optimisation analysis 
  begins with the identification of some natural candidates $\bu$ to the role of optimal controls.
  Secondly, for such controls $\bu$ the forward and the backward systems are solved, 
  so that the respective variables $\varphi=\varphi(\bu)$ and $\nabla P=\nabla P(\bu)$
  are identified. Finally, if condition \eqref{NC2} is not met, then the candidate $\bu$
  is cut off from the analysis, otherwise it is confirmed.
  Nonetheless, let stress again that condition \eqref{NC2} is only a necessary requirement,
  and can only help to restrict the class of potential optimal controls. 
  In order to further refine the analysis, sufficient conditions for optimality should be investigated.
  The mathematical idea behind this is very natural: if the reduced cost functional $\widetilde J$
  can be shown to be twice (Fr\'echet or G\^ateaux) differentiable, then 
  any control $\bu$ satisfying the first order--stationary condition \eqref{NC2} and 
  the positive definiteness condition $D^2\widetilde J(\bu)>0$ is an optimal control.
  Such second--order analysis is extremely challenging, and to the best of the 
  author's knowledge it has been performed so far only in relation to some selected 
  optimal control problems in the deterministic setting 
  \cite{col-far-hass-gil-spr2, colli-sprek-optACDBC}. 
  In the stochastic case, the second--order analysis is open, and is currently 
  being investigated in a work in preparation. 
\end{remark}

\section{Well-posedness of the state system}
\label{sec:WP}
This section is devoted to
the proof of Theorem~\ref{thm1} about well-posedness of the state system.

\subsection{Uniqueness}
\label{ssec:uniq}
Let $\{\bu_i\}_{i=1,2}\subset\mathcal U$ and let us denote by
$\{(\varphi_i, \mu_i)\}_{i=1,2}$
any respective solutions to \eqref{eq1}--\eqref{eq4} in the sense of Theorem~\ref{thm1}.
Let us set for brevity of notation $\varphi:=\varphi_1-\varphi_2$, $\mu:=\mu_1-\mu_2$,
$\bu:=\bu_1-\bu_2$: then we have
\[
  \d\varphi - \Delta \mu\,\d t + \bu\cdot\nabla\varphi_1\,\d t 
  + \bu_2\cdot\nabla\varphi\,\d t = (B(\varphi_1)-B(\varphi_2))\,\d W\,,
  \qquad\varphi(0)=0\,,
\]
where the equality is intended in the usual variational sense
of Theorem~\ref{thm1}.

Taking $\frac1{|\OO|}\in V_1$ as test function yields directly
by assumption {\bf A3} that 
$\varphi_\OO=0$, so that actually $\varphi\in L^p_\cP(\Omega; C^0([0,T]; V_{1,0}^*))$
and $B(\varphi_1)-B(\varphi_2)\in L^p_\cP(\Omega; L^2(0,T; \cL^2(K,V_{1,0}^*)))$.
Hence, It\^o's formula for the function $\frac12\norm{\nabla\mN\varphi}_H^2$
yields
\begin{align*}
  &\frac12\norm{\nabla\mathcal N\varphi(t)}_{H}^2
  +\int_{Q_t}|\nabla\varphi|^2 + \int_{Q_t}(\Psi'(\varphi_1)-\Psi'(\varphi_2))\varphi
  +\int_{Q_t}(\bu\cdot\nabla\varphi_1+ \bu_2\cdot\nabla\varphi)
  \mathcal N\varphi\\
  &=\frac12\int_0^t\norm{\nabla\mathcal N (B(\varphi_1)-B(\varphi_2))(s)}_{\cL^2(K,H)}^2\,\d s
  +\int_0^t\left(\mathcal N\varphi(s),
  (B(\varphi_1)-B(\varphi_2))(s)\,\d W(s)\right)_H\,.
\end{align*}
Now, the mean value theorem and assumption {\bf A1} give
\[
  \int_{Q_t}(\Psi'(\varphi_1)-\Psi'(\varphi_2))\varphi 
  \geq - C_\Psi\int_{Q_t}|\varphi|^2\,,
\]
while the inclusion $V_1\embed L^6(\OO)$,
the H\"older and the Poincar\'e-Wirtinger inequalities yield
\begin{align*}
  &\int_{Q_t}(\bu\cdot\nabla\varphi_1+ \bu_2\cdot\nabla\varphi)\mathcal N\varphi\leq
  c\int_0^t\left(\norm{\nabla\varphi_1(s)}_H\norm{\bu(s)}_U+
  \norm{\bu_2(s)}_U\norm{\nabla\varphi(s)}_H\right)
  \norm{\mathcal N\varphi(s)}_{V_1}\,\d s\\
  &\leq 
  \norm{\varphi_1}_{L^\infty(0,T; V_1)}^2\norm{\bu}_{L^2(0,T; U)}^2
  +\frac12\int_{Q_t}|\nabla\varphi|^2
  +c\int_0^t\left(1+ \norm{\bu_2(s)}_U^2\right)\norm{\nabla\mathcal N\varphi(s)}^2_{H}\,\d s\,.
\end{align*}
Furthermore, assumption {\bf A3} ensure that 
\[
  \int_0^t\norm{\nabla\mathcal N (B(\varphi_1)-B(\varphi_2))(s)}_{\cL^2(K,H)}^2\,\d s
  \leq c\int_{Q_t}|\varphi|^2\,.
\]
Using the compactness inequality \eqref{comp_ineq} and rearranging 
the terms we are left with 
\begin{align}
  \nonumber
  \norm{\nabla\mathcal N\varphi(t)}_{H}^2
  +\int_{Q_t}|\nabla\varphi|^2
  &\leq c\norm{\varphi_1}_{L^\infty(0,T; V_1)}^2\norm{\bu}_{L^2(0,T; U)}^2+
  c\int_0^t\left(1+ \norm{\bu_2(s)}_U^2\right)\norm{\nabla\mathcal N\varphi(s)}^2_{H}\,\d s\\
  \label{uniq_aux}
  &+c\int_0^t\left(\mathcal N\varphi(s),
  (B(\varphi_1)-B(\varphi_2))(s)\,\d W(s)\right)_H\,.
\end{align}
On the right--hand side we have, by the H\"older inequality in time, 
\[
  \int_0^t\left(1+ \norm{\bu_2(s)}_U^2\right)\norm{\nabla\mathcal N\varphi(s)}^2_{H}\,\d s
  \leq c t^{1-\frac2p}\left(1+\norm{\bu}_{L^p(0,T; U)}^2\right)
  \norm{\nabla\mathcal N\varphi}^2_{L^\infty(0,t;H)}\,,
\]
and, thanks to the Burkh\"older-Davis-Gundy and
the Young inequalities, assumption {\bf A3},
and again the compactness inequality \eqref{comp_ineq},
\begin{align*}
  \E\sup_{r\in[0,t]}\left|\int_0^r\left(\mathcal N\varphi(s),
  (B(\varphi_1)-B(\varphi_2))(s)\,\d W(s)\right)_H\right|^{p/2}
  \leq \frac18\E\norm{\nabla\mathcal N\varphi}_{L^\infty(0,t; H)}^p
  + c \E\norm{\varphi}^p_{L^2(0,t; H)}&\\
  \leq
  \frac18\E\norm{\nabla\mathcal N\varphi}_{L^\infty(0,t; H)}^p+
  \frac12\E\norm{\nabla\varphi}^p_{L^2(0,t; H)}
  + c \E\norm{\nabla\mN\varphi}^p_{L^2(0,t; H)}&\,.
\end{align*}
Consequently,
taking power $p/2$ at both sides of \eqref{uniq_aux}
and rearranging the terms yield
\begin{align*}
  &\E\norm{\nabla\mathcal N\varphi}_{L^\infty(0,t;H)}^p
  +\E\norm{\nabla\varphi}_{L^2(0,t; H)}^p\\
  &\leq c\norm{\bu}_{\mU}^p\E\norm{\varphi_1}_{L^\infty(0,T; V_1)}^p + 
  c t^{\frac{p}2-1}(1+\norm{\bu_2}_\mU^p)\E\norm{\nabla\mathcal N\varphi}_{L^\infty(0,t;H)}^p\,.
\end{align*}
Hence, setting 
\[
T_0:=\left(\frac12c^{-1}(1+\norm{\bu_2}_\mU^p)^{-1}\right)^{\frac2{p-2}}\wedge T\,,
\]
we get
\[
  \E\norm{\nabla\mathcal N\varphi}_{L^\infty(0,T_0;H)}^p
  +\E\norm{\nabla\varphi}_{L^2(0,T_0; H)}^p\\
  \leq c\norm{\bu}_{\mU}^p\E\norm{\varphi_1}_{L^\infty(0,T; V_1)}^p + 
  \frac12\E\norm{\nabla\mathcal N\varphi}_{L^\infty(0,T_0;H)}^p\,.
\]
Since $T_0$ is independent of the initial time, 
we can iterate the procedure and close the estimate 
on each subinterval $[kT_0, (k+1)T_0]$ for all $k\in\enne$
until $(k+1)T_0>T$: 
summing up, noting that the number of such subintervals is
less than $\frac{T}{T_0}+1$,
and renominating $c$ independently of $\bu_2$,
we get then
\begin{align*}
  &\norm{\varphi_1-\varphi_2}_{L^p_\cP(\Omega; C^0([0,T]; V_1^*)\cap L^2(0,T; V_1))}^p
  \leq c\norm{\varphi_1}_{L^p(\Omega;L^\infty(0,T; V_1))}^p
  \left(1 + \norm{\bu_2}_\mU^{\frac{2p}{p-2}}\right)
  \norm{\bu_1-\bu_2}^p_{\mU}\,,
\end{align*}
from which uniqueness of solutions follows.

\subsection{Approximation}
\label{ssec:approx}
We turn now to existence of solutions.
First of all, for every $\lambda$ let $\beta_\lambda:\erre\to\erre$ be the Yosida approximation 
of $\beta$ and $\widehat\beta_\lambda:\erre\to[0,+\infty)$ be the 
Moreau-Yosida regularisation of $\widehat\beta$, which are
defined, respectively, as 
\[
  \beta_\lambda(r):=\frac{r-(I+\lambda\beta)^{-1}(r)}{\lambda}\,,
  \qquad
  \widehat\beta_\lambda(r):=\int_0^r\beta_\lambda(s)\,\d s\,, \qquad r\in\erre.
\]
Let us recall that $\beta_\lambda$ is $\frac1\lambda$--Lipschitz continuous,
$\widehat\beta_\lambda$ is convex and quadratic at $\infty$, and 
as $\lambda\searrow0$ it holds that $\beta_\lambda(r)\to\beta(r)$ and 
$\widehat\beta_\lambda(r)\nearrow\widehat\beta(r)$ for all $r\in\erre$.
For further details about the properties of $\beta_\lambda$ and $\widehat\beta_\lambda$
we refer to the monograph \cite[Ch.~2]{barbu-monot}.
We define the approximated double-well potential as
\[
  \Psi_\lambda:\erre\to\erre\,, \qquad 
  \Psi_\lambda(r):=\Psi(0) + \widehat\beta_\lambda(r) - \frac{C_\Psi}{2}r^2\,, \quad r\in\erre\,,
\]
so that in particular we have $\Psi_\lambda'(r)=\beta_\lambda(r) - C_\Psi r$ for $r\in\erre$.
Secondly, we define
\[
  \bu_\lambda:= \rho_\lambda*\bu\,,
\]
where $(\rho_\lambda)_\lambda\subset C^\infty_c(\erre)$ is
is a classical non--anticipative 
sequence of mollifiers in time. In particular, let us point out that it holds 
\[
  \bu_\lambda\in L^\infty_\cP(\Omega\times(0,T); U)\,,
  \qquad\bu_\lambda\to\bu \quad\text{in }L^q_\cP(\Omega; L^p(0,T; U))\quad\forall\,q\geq1\,.
\]

The approximated system is obtained by replacing $\Psi'$ with $\Psi'_\lambda$
and $\bu$ with $\bu_\lambda$
in \eqref{eq1}--\eqref{eq4}:
\begin{align}
  \label{eq1_lam}
  \d \varphi_\lambda - \Delta\mu_\lambda\,\d t + \bu_\lambda\cdot\nabla\varphi_\lambda\,\d t =
  B(\varphi_\lambda)\,\d W \qquad&\text{in } (0,T)\times\OO\,,\\
  \label{eq2_lam}
  \mu_\lambda = -\Delta\varphi_\lambda + 
  \Psi_\lambda'(\varphi_\lambda) \qquad&\text{in } (0,T)\times\OO\,,\\
  \label{eq3_lam}
  {\bf n}\cdot\nabla\varphi_\lambda = {\bf n}\cdot\nabla\mu_\lambda= 0 
  \qquad&\text{in } (0,T)\times\partial\OO\,,\\
  \label{eq4_lam}
  \varphi_\lambda(0)=\varphi_0 \qquad&\text{in } \OO\,.
\end{align}
We formulate \eqref{eq1_lam}--\eqref{eq4_lam} in an abstract way as
\beq\label{eq_ab_lam}
  \d \varphi_\lambda +
  (\mathcal A_\lambda + \mathcal C_\lambda)(\varphi_\lambda)\,\d t =
  B(\varphi_\lambda)\,\d W\,, \qquad
  \varphi_\lambda(0)=\varphi_0\,,
\eeq
where the variational operators
\[
  \mathcal A_\lambda:V_2\to V_2^*\,,\qquad
  \mathcal C_\lambda:\Omega\times[0,T]\times V_2\to V_2^*\,,
\]
are defined as
\[
  \ip{\mathcal A_\lambda(y)}{\zeta}
  :=\int_\OO(-\Delta \zeta)(-\Delta y + \Psi'_\lambda(y))\,, \quad y,\zeta\in V_2\,,\\
\]
and
\[
  \ip{\mathcal C_\lambda(\omega,t,y)}{\zeta}:= - \int_\OO y \bu_\lambda(\omega,t)
  \cdot\nabla\zeta\,, \qquad y,\zeta\in V_2\,,\quad t\in[0,T]\,.
\]
Since $\Psi'_\lambda$ is Lipschitz-continuous, it is not difficult 
to show (see for example \cite[Lem.~3.1]{scar-SCH}) that 
$\mathcal A_\lambda$ is weakly monotone, 
weakly coercive, and linearly bounded, in the sense that there are
two constants $c_\lambda, c_\lambda'>0$ such that 
\[
  \ip{\mathcal A_\lambda(y_1)-\mathcal A_2(y_2)}{y_1-y_2}
  \geq c_\lambda\norm{y_1-y_2}_{V_2}^2 - c_\lambda'\norm{y_1-y_2}_H^2
  \qquad\forall\,y_1,y_2\in V_2
\]
and 
\[
  \norm{\mathcal A_\lambda(y)}_{V_2^*} \leq c_\lambda'(1+\norm{y}_{V_2})
  \qquad\forall\,y\in V_2\,.
\]
As far as the convection 
operator $\mathcal C_\lambda$ is concerned, 
since $\operatorname{div}\bu_\lambda=0$, 
thanks to the divergence theorem we have
\begin{align*}
  \ip{\mathcal C_\lambda(y_1)-\mathcal C_\lambda(y_2)}{y_1-y_2}&=
  - \int_\OO (y_1-y_2) \bu_\lambda\cdot\nabla(y_1-y_2)=0\,,
\end{align*}
and, thanks to the H\"older inequality and the inclusion $V_1\embed L^6(\OO)$,
\[
  \norm{C_\lambda(y)}_{V_2^*}=
  \sup_{\norm{\zeta}_{V_2}\leq 1}\left\{- \int_\OO y \bu_\lambda\cdot\nabla\zeta\right\}
  \leq\norm{y}_{H}\norm{\bu_\lambda}_{U}\leq 
  \norm{\bu_\lambda}_{L^\infty_\cP(\Omega\times(0,T); U)}\norm{y}_{V_2}
  \quad\forall\,y\in V_2\,.
\]
Hence, the operator 
$\mathcal A_\lambda+\mathcal C_\lambda:\Omega\times[0,T]\times V_2\to V_2^*$
is weakly monotone, weakly coercive, and linearly bounded.
Besides, due to the Lipschitz-continuity of $\Psi_\lambda'$
and the regularity of $\bu_\lambda$
it is immediate to check that it is also hemicontinuous.
Moreover, assumption {\bf A3} ensure that $B: H\to\cL^2(K,H)$
is Lipschitz-continuous.
It follows then by the classical variational approach to SPDEs
by Pardoux \cite{Pard} and Krylov--Rozovskii \cite{KR-SPDE} that 
the evolution equation \eqref{eq_ab_lam} admits a unique variational solution
\[
  \varphi_\lambda \in L^2_\cP(\Omega; C^0([0,T]; H)\cap L^2(0,T; V_2))\,.
\]
Let us set $\mu_\lambda:=-\Delta\varphi_\lambda + \Psi_\lambda'(\varphi_\lambda)$
as the approximated chemical potential.

\subsection{Uniform estimates}
\label{ssec:est}
It\^o's formula for the square of the $H$-norm yields 
\begin{align*}
  &\frac12\norm{\varphi_\lambda(t)}_H^2
  +\int_{Q_t}|\Delta\varphi_\lambda|^2 
  + \int_{Q_t}\Psi_\lambda'(\varphi_\lambda)(-\Delta\varphi_\lambda)
  -\int_{Q_t}\varphi_\lambda \bu_\lambda\cdot\nabla\varphi_\lambda\\
  &=\frac12\norm{\varphi_0}_H^2
  +\frac12\int_0^t\norm{B(\varphi_\lambda(s))}_{\cL^2(K,H)}^2\,\d s
  +\int_0^t\left(\varphi_\lambda(s), B(\varphi_\lambda(s))\,\d W(s)\right)_H\,.
\end{align*}
Now, on the left-hand side we have, thanks to the monotonicity of $\beta_\lambda$,
\[
   \int_{Q_t}\Psi_\lambda'(\varphi_\lambda)(-\Delta\varphi_\lambda)=
   \int_{Q_t}\beta_\lambda'(\varphi_\lambda)|\nabla\varphi_\lambda|^2
   -C_\Psi\int_{Q_t}\varphi_\lambda(-\Delta\varphi_\lambda)
   \geq -\frac14\int_{Q_t}|\Delta\varphi_\lambda^2| - C_\Psi^2\int_{Q_t}|\varphi_\lambda|^2\,.
\]
Also, by the H\"older inequality and the inclusion $V_1\embed L^6(\OO)$ it holds
\begin{align*}
  -\int_{Q_t}\varphi_\lambda \bu_\lambda\cdot\nabla\varphi_\lambda
  &\geq -\int_0^t\norm{\varphi_\lambda(s)}_H\norm{\bu_\lambda(s)}_{U}
  \norm{\varphi(s)}_{V_2}\,\d s\,.
\end{align*}
Thanks to the elliptic regularity theory for the Neumann problem 
(see e.g.~\cite[\S 9.6]{brezis-anal})
there is $c>0$ independent of $\lambda$ such that $\norm{\zeta}_{V_2}\leq
c(\norm{\zeta}_H+\norm{\Delta\zeta}_H)$ for every $\zeta\in V_2$: consequently, 
renominating $c$ and using the Young inequality we get
\[
  -\int_{Q_t}\varphi_\lambda \bu_\lambda\cdot\nabla\varphi_\lambda
  \geq -\frac14\int_{Q_t}|\Delta\varphi_\lambda|^2 
  -c^2\int_0^t\norm{\varphi_\lambda(s)}_H^2(1+\norm{\bu_\lambda(s)}_{U}^2)\,\d s\,.
\]
Furthermore, noting that $\frac{2\gamma}{\gamma-1}\geq4$
since $\gamma\in[1,2]$, assumption {\bf A3} yields
\[
  \frac12\int_0^t\norm{B(\varphi_\lambda(s))}_{\cL^2(K,H)}^2\,\d s
  \leq c\,.
\]
Putting this information together and using 
assumption  on the right-hand side we get,
possibly updating the value of $c$,
\begin{align*}
  \frac12\norm{\varphi_\lambda(t)}_H^2
  +\frac12\int_{Q_t}|\Delta\varphi_\lambda|^2
  &\leq\frac12\norm{\varphi_0}_H^2 
  +c\int_0^t\norm{\varphi_\lambda(s)}_H^2(1+\norm{\bu_\lambda(s)}_{U}^2)\,\d s\\
  &+\int_0^t\left(\varphi_\lambda(s), B(\varphi_\lambda(s))\,\d W(s)\right)_H
  \qquad\forall\,t\in[0,T]\,,\quad\P\text{-a.s.}
\end{align*}
Taking now power $p/2$ at both sides, 
the stochastic integral on the right-hand side can be treated
again thanks to {\bf A3},
using classical computations based on the Burkholder-Davis-Gundy inequality 
(see for example \cite[Lem.~4.3]{mar-scar-diss}). Consequently, 
the same iterative argument used in Subsection~\ref{ssec:uniq}
ensures that 
\beq
  \label{est1}
  \norm{\varphi_\lambda}_{L^p_\cP(\Omega; C^0([0,T]; H)\cap L^2(0,T; V_2))}^p \leq 
  c\left(1+\norm{\bu}_{\mU}^{\frac{2p}{p-2}}\right)\,.
\eeq

In order to deduce further estimates on $\varphi_\lambda$ and $\mu_\lambda$,
we rely on the free-energy estimate. Namely, we consider the 
approximated energy 
\[
  \zeta\mapsto E_\lambda(\zeta):=
  \frac12\int_\OO|\nabla\zeta|^2 + \int_\OO\Psi_\lambda(\zeta)\,,
  \qquad \zeta\in V_1\,.
\]
Clearly, $E_\lambda$ is well-defined and of class $C^1$ in $V_1$, with derivative
\[
  DE_\lambda:V_1\to V_1^*\,, \qquad 
  DE_\lambda(\zeta)=\mathcal L\zeta + \Psi'_\lambda(\zeta)\,, \quad\zeta\in V_1\,,
\]
so that in particular we have $DE_\lambda(\varphi_\lambda)=
\mu_\lambda$.
Moreover, the Lipschitz-continuity of $\Psi'_\lambda$ ensures that 
$DE_\lambda:V_1\to V_1^*$ is actually Fr\'echet-differentiable with 
\[
  D^2E_\lambda(\zeta)[z_1,z_2]=\int_\OO\nabla z_1\cdot\nabla z_2
  +\int_\OO\Psi_\lambda''(\zeta)z_1z_2\,, \quad\zeta,z_2,z_2\in V_1\,.
\]
Now, we would like to write It\^o's formula for $E_\lambda(\varphi_\lambda)$:
in order to do this, we need to show first that $\varphi_\lambda$ and $\mu_\lambda$
enjoy more regularity. This can be shown 
by performing a further approximation on the problem 
(for example, the classical 
Faedo--Galerkin approximation of the abstract evolution equation 
\eqref{eq_ab_lam}). Indeed, 
by the classical variational theory on stochastic evolution equations \cite{LiuRo},
there is a sequence $(H_n)_n$ of finite--dimensional subspaces of $H$,
included in $V_2$ and with $\cup_nH_n$ dense in $H$,
such that, setting $P_n:V_2^*\to H_n$ as the orthogonal projection onto $H_n$,
the unique solution $(\varphi_\lambda^n, \mu_\lambda^n)$ of the finite dimensional system
\begin{align*}
  \d \varphi_\lambda^n - \Delta\mu^n_\lambda\,\d t + 
  P_n(\bu_\lambda\cdot\nabla\varphi_\lambda^n)\,\d t =
  P_nB(\varphi_\lambda^n)\,\d W \qquad&\text{in } (0,T)\times\OO\,,\\
  \mu_\lambda^n = -\Delta\varphi_\lambda^n + 
  P_n\Psi_\lambda'(\varphi_\lambda^n) \qquad&\text{in } (0,T)\times\OO\,,\\
  {\bf n}\cdot\nabla\varphi_\lambda^n = {\bf n}\cdot\nabla\mu_\lambda^n= 0 
  \qquad&\text{in } (0,T)\times\partial\OO\,,\\
  \varphi^n_\lambda(0)=\varphi_0^n \qquad&\text{in } \OO\,,
\end{align*}
satisfy, as $n\to\infty$,
\[
  \varphi_\lambda^n\wto\varphi_\lambda \quad\text{in } L^p_\cP(\Omega; L^2(0,T; V_2))\,,
  \qquad
  \mu_\lambda^n\wto\mu_\lambda \quad\text{in } L^p_\cP(\Omega; L^2(0,T; H))\,.
\]
At this point, the finite--dimensional It\^o formula for ${E_\lambda}_{|H_n}$ yields 
\begin{align*}
  &\frac12\int_\OO|\nabla\varphi_\lambda^n(t)|^2 
  +\int_\OO\Psi_\lambda(\varphi_\lambda^n(t))
  +\int_{Q_t}|\nabla\mu_\lambda^n|^2 
  =
  \frac12\int_\OO|\nabla\varphi_0^n|^2 
  +\int_\OO\Psi_\lambda(\varphi_0^n)
  +\int_{Q_t}\varphi_\lambda^n \bu_\lambda\cdot\nabla\mu_\lambda^n\\
  &+\frac12\int_0^t
  \norm{\nabla P_nB(\varphi_\lambda^n(s))}^2_{\cL^2(K,H)}\,\d s
  +\sum_{j=0}^\infty\int_{Q_t}
  \Psi_\lambda''(\varphi_\lambda^n)|P_nB(\varphi_\lambda^n)e_j|^2
  +\int_0^t\left(\mu_\lambda^n(s), B(\varphi_\lambda^n(s))\,\d W(s)\right)_H
\end{align*}
for every $t\in[0,T]$, $\P$-almost surely.
We show now uniform estimates on the terms on the right--hand side,
independent of both $\lambda$ and $n$. These will show a posteriori 
that $\varphi_\lambda$ and $\mu_\lambda$ are actually more regular. 
For this reason and for brevity of notation, we omit from now on 
the dependence on $n$ and refer to \cite{scar-SCH, scar-SVCH} for  more detail.

To this end, noting that the definition of $\mu_\lambda$ and assumption {\bf A1} imply
\[
  |(\mu_\lambda)_\OO| =|(\Psi_\lambda'(\varphi_\lambda))_\OO| \leq
  \norm{\Psi_\lambda'(\varphi_\lambda)}_{L^1(\OO)}
  \leq c\left(1 + \int_\OO\Psi_\lambda(\varphi_\lambda)\right)\,, 
\]
on the left-hand side we get 
\[
  \int_\OO\Psi_\lambda(\varphi_\lambda(t)) \geq \frac1c|(\mu_\lambda(t))_\OO| - c\,.
\]
On the right-hand side, thanks to the 
H\"older and Young inequalities, the inclusion $V_1\embed L^6(\OO)$,
and the estimate \eqref{est1},
proceeding as in Subsection~\ref{ssec:uniq} we have
\begin{align*}
  \int_\OO\Psi_\lambda(\varphi_0)
  +\int_{Q_t}\varphi_\lambda \bu_\lambda\cdot\nabla\mu_\lambda&\leq
  \int_\OO\Psi(\varphi_0) + \frac12\int_{Q_t}|\nabla\mu_\lambda|^2
  +\frac12\int_0^t\norm{\varphi_\lambda(s)}^2_{V_1}\norm{\bu_\lambda(s)}_U^2\,\d s\\
  &\leq c+\frac12\int_{Q_t}|\nabla\mu_\lambda|^2
  +ct^{1-\frac{2}{p}}\norm{\bu}_\mU^2\norm{\nabla\varphi}_{L^\infty(0,t; H)}^2
\end{align*}
Moreover, assumptions {\bf A3} and {\bf A1} yield,
together with the H\"older inequality and \eqref{est1},
\begin{align*}
  &\frac12\int_0^t
  \norm{\nabla B(\varphi_\lambda(s))}^2_{\cL^2(K,H)}\,\d s
  +\sum_{j=0}^\infty\int_{Q_t}
  \Psi_\lambda''(\varphi_\lambda)|B(\varphi_\lambda)e_j|^2\\
  &
  \leq c\left(1+\int_0^t\norm{\varphi_\lambda(s)}_{V_1}^2\,\d s\right)
  +\sum_{j=0}^\infty
  \int_0^t\norm{\Psi_\lambda''(\varphi_\lambda(s))}_{L^\gamma(\OO)}
  \norm{B(\varphi_\lambda(s))e_j}_{L^{\frac{2\gamma}{\gamma-1}}(\OO)}^2\,\d s\\
  &\leq c\left(1 + \int_0^t\norm{\nabla\varphi_\lambda(s)}_{H}^2\,\d s
  +\int_{Q_t}\Psi_\lambda(\varphi_\lambda)\right)\\
  &\leq c + ct\norm{\nabla\varphi_\lambda}_{L^\infty(0,t; H)}^2+
  ct\norm{\Psi_\lambda(\varphi_\lambda)}_{L^\infty(0,t; L^1(\OO))}\,.
\end{align*}
Finally, the Burkholder--Davis--Gundy and 
the Poincar\'e--Wirtinger inequalities give, together with assumption {\bf A3},
\begin{align*}
  \E\sup_{r\in[0,t]}\left|\int_0^r\left(\mu_\lambda(s), 
  B(\varphi_\lambda(s))\,\d W(s)\right)_H\right|^{p/2}
  &\leq c\E\left(\int_0^t\norm{\mu_\lambda(s)}_H^2
  \norm{B(\varphi_\lambda(s))_H}_{\cL^2(K,H)}^2\,\d s\right)^{p/4}\\
  \leq c\E\norm{\mu_\lambda}_{L^2(0,t; H)}^{p/2}
  &\leq \delta\E\norm{\nabla\mu_\lambda}^p_{L^2(0,t; H)}
  + c_\delta\left(1 + \E\norm{(\mu_\lambda)_\OO}^{p/2}_{L^2(0,t)}\right)\,,
\end{align*}
for every $\delta>0$,
where we have updated the value of $c$ and $c_\delta$ step--by--step, 
independently of $\lambda$.
Putting all this information together, 
choosing $\delta$ sufficiently small, rearranging the terms,
and updating again the value of $c$, we infer that
\begin{align*}
  &\E\norm{\nabla\varphi_\lambda}_{L^\infty(0,t; H)}^p 
  +\E\norm{\Psi_\lambda(\varphi_\lambda)}_{L^\infty(0,t;L^1(\OO))}^{p/2}
  +\E\norm{(\mu_\lambda)_\OO}_{L^\infty(0,t)}^{p/2}
  +\E\norm{\nabla\mu_\lambda}^{p}_{L^2(0,t; H)} \\
  &\leq c\left[1 + 
  \left(t^{\frac{p}2-1}\norm{\bu}_\mU^p+t^{\frac{p}2}\right)
  \E\norm{\nabla\varphi}_{L^\infty(0,t; H)}^p\right.\\
  &\qquad\left.+
  t^{\frac{p}{2}}\E\norm{\Psi_\lambda(\varphi_\lambda)}_{L^\infty(0,t; L^1(\OO))}^{p/2}
  +t^{p/4}\E\norm{(\mu_\lambda)_\OO}_{L^\infty(0,t)}^{p/2}\right]
  \qquad\forall\,t\in[0,T]\,.
\end{align*}
Consequently, we can close the estimate on
a certain subinterval $[0,T_0]$, where $T_0$ is chosen sufficiently small
in order to incorporate the terms on the right-hand side into the 
corresponding ones on the left. Also, 
a patching argument as in Subsection~\ref{ssec:uniq} allows then to 
extend the estimate to the whole interval $[0,T]$, and we obtain
\begin{align}
  \nonumber
  \norm{\varphi_\lambda}_{L^p(\Omega; L^\infty(0,T; V_1))}
  &+\norm{\mu_\lambda}_{L^{p/2}_\cP(\Omega; L^2(0,T; V_1))}
  +\norm{\nabla\mu_\lambda}_{L^{p}_\cP(\Omega; L^2(0,T; H))}\\
  \label{est2}
  &+\norm{\Psi_\lambda(\varphi_\lambda)}_{L^{p/2}(\Omega; L^\infty(0,T; L^1(\Omega)))} 
  \leq c\left(1+\norm{\bu}_{\mU}^{\frac{2p}{p-2}}\right)\,,
\end{align}
which by comparison in $\mu_\lambda=-\Delta\varphi_\lambda + \Psi_\lambda'(\varphi_\lambda)$
and estimate \eqref{est1} gives also
\beq
  \label{est3}
  \norm{\Psi'_\lambda(\varphi_\lambda)}_{L^{p/2}_\cP(\Omega; L^2(0,T; H))} \leq 
  c\left(1+\norm{\bu}_{\mU}^{\frac{2p}{p-2}}\right)\,.
\eeq

Finally, note that by assumption {\bf A3} and the estimate \eqref{est2} we have 
\[
  \norm{B(\varphi_\lambda)}_{L^\infty(\Omega\times(0,T); \cL^2(K,H))
  \cap L^p(\Omega; L^\infty(0,T; \cL^2(K,V_1)))}\leq c\,,
\]
so that the classical result by Flandoli \& Gatarek \cite[Lem.~2.1]{fland-gat} ensures 
in particular that 
\beq
  \label{est4}
  \norm{I_\lambda:=
  \int_0^\cdot B(\varphi_\lambda(s))\,\d W(s)}_{L^p_\cP(\Omega; W^{s,p}(0,T; V_1))} 
  \leq c_{s}
  \qquad\forall\,s\in(0,1/2)\,.
\eeq
Consequently, by comparison in \eqref{eq1_lam} it is not difficult to check that 
\[
  \norm{\varphi_\lambda}_{L^p_\cP(\Omega; 
  W^{1,2}(0,T; V_1^*) +W^{s,p}(0,T; V_1))} \leq c_s
  \qquad\forall\,s\in(0,1/2)\,.
\]
Now, recalling that $p>2$, for all arbitrary 
$s\in (0,1/2)$ we have that 
$s-\frac1p\leq\frac12$,
so that 
the usual Sobolev embeddings ensure that
\[
  W^{1,2}(0,T; V_1^*)\embed W^{s,p}(0,T; V_1^*) \quad\forall\,s\in(0,1/2)\,,
\]
and we deduce that 
\beq
  \label{est5}
  \norm{\varphi_\lambda}_{L^p_\cP(\Omega; W^{s,p}(0,T; V_1^*))} \leq c_s
  \qquad\forall\,s\in(0,1/2)\,.
\eeq

\subsection{Passage to the limit}
\label{ssec:pass_lim}
From the estimates \eqref{est1}--\eqref{est3}, 
there exists a pair $(\varphi,\mu)$, with 
\[
  \varphi\in L^p_w(\Omega; L^\infty(0,T; V_1))\cap L^p_\cP(\Omega; L^2(0,T; V_2))\,,
  \qquad \mu\in L^{p/2}_\cP(\Omega; L^2(0,T; V_1))
\]
such that, as $\lambda\searrow0$, on a non-relabelled subsequence we have
\begin{align*}
  \varphi_\lambda \wstarto \varphi \qquad&\text{in } L^p_w(\Omega; L^\infty(0,T; V_1))
  \cap L^p_\cP(\Omega; L^2(0,T; V_2))\,,\\
  \mu_\lambda \wto \mu \qquad&\text{in } L^{p/2}_\cP(\Omega; L^2(0,T; V_1))\,.
\end{align*}
Now, since $p>2$, we can fix $\bar s\in (\frac1p, \frac12)$, so that $\bar s p>1$:
with this choice, 
by the classical Aubin--Lions--Simon compactness results \cite[Cor.~5]{simon} we have
\[  
  L^\infty(0,T; V_1)\cap L^2(0,T; V_2)\cap W^{\bar s,p}(0,T; V_1^*)\embed
  C^0([0,T]; H)\cap L^2(0,T; V_1) \qquad\text{compactly}\,.
\]
Hence, setting $\mathcal B_n$ as the closed ball of radius $n$ in 
$L^\infty(0,T; V_1)\cap L^2(0,T; V_2)\cap W^{\bar s,p}(0,T; V_1^*)$,
we have that $\mathcal B_n$ 
is compact in $C^0([0,T]; H)\cap L^2(0,T; V_1)$, for every $n\in\enne$.
Consequently, denoting by $\nu_\lambda$ the law of $\varphi_\lambda$
on $C^0([0,T]; H)\cap L^2(0,T; V_1)$ for brevity, 
the Markov inequality and the uniform estimates \eqref{est1}, \eqref{est2}, and \eqref{est5}
yield 
\begin{align*}
  \nu_\lambda(\mathcal B_n^c)&=
  \P\{\norm{\varphi_\lambda}_{L^\infty(0,T; V_1)
  \cap L^2(0,T; V_2)\cap W^{\bar s,p}(0,T; V_1^*)}>n\}\\
  &\leq\frac1n\E\norm{\varphi_\lambda}_{L^\infty(0,T; V_1)
  \cap L^2(0,T; V_2)\cap W^{\bar s,p}(0,T; V_1^*)}\leq\frac{c}{n}\,,
\end{align*}
from which 
\[
  \lim_{n\to\infty}\sup_{\lambda>0}\nu_\lambda(\mathcal B_n^c)=0\,.
\]
By the Prokhorov theorem this implies
that 
\[ 
  \text{the laws of } (\varphi_\lambda)_\lambda \text{ are tight on }
  C^0([0,T]; H)\cap L^2(0,T; V_1)\,.
\]
Similarly, estimate \eqref{est4} ensures by the same argument that 
\[
  \text{the laws of } (I_\lambda)_\lambda \text{ are tight on }
  C^0([0,T]; H)\,.
\]
Let us show now that, possibly on a further subsequence, 
we have also the strong convergence
\beq
  \label{strong}
  \varphi_\lambda\to\varphi \qquad\text{in } C^0([0,T]; H)\cap L^2(0,T; V_1)
  \quad\P\text{-a.s.}
\eeq
To this end, we use the following lemma due to 
Gy\"ongy \& Krylov \cite[Lem.~1.1]{gyo-kry},
which characterises the convergence in probability in a Polish space.
\begin{lem}\label{lem:gyo-kry}
  Let $\mathcal X$ be a Polish space and $(Z_n)_n$ be a sequence
  of $\mathcal X$-valued random variables. Then $(Z_n)_n$ converges
  in probability if and only if
  for any pair of subsequences $(Z_{n_k})_k$ and $(Z_{n_j})_j$, there exists 
  a joint sub-subsequence $(Z_{n_{k_\ell}}, Z_{n_{j_\ell}})_\ell$ converging 
  in law to a probability measure $\nu$ on $\mathcal X\times\mathcal X$ such
  that $\nu(\{(z_1,z_2)\in\mathcal X\times\mathcal X: z_1=z_2\})=1$.
\end{lem}
We apply this lemma to $\mathcal X=C^0([0,T]; H)\cap L^2(0,T; V_1)$
and $(\varphi_{\lambda})_\lambda$. Given two arbitrary subsequences
$(\varphi_{\lambda_k})_k$ and $(\varphi_{\lambda_j})_j$, since 
the laws of the pairs $(\varphi_{\lambda_k}, \varphi_{\lambda_j})_{k,j}$ are tight
on $(C^0([0,T]; H)\cap L^2(0,T; V_1))^2$, there is a joint 
subsequence $(\varphi_{\lambda_{k_i}}, \varphi_{\lambda_{j_i}})_i$
converging weakly to a probability measure $\nu$ on
$(C^0([0,T]; H)\cap L^2(0,T; V_1))^2$.
By the Skorokhod representation theorem 
\cite[Thm.~2.7]{ike-wata} 
and \cite[Thm.~1.10.4, Add.~1.10.5]{vaa-well}, there exist
a new probability space $(\Omega', \cF', \P')$ and 
measurable maps $\phi_i:(\Omega', \cF')\to(\Omega,\cF)$,
such that $\P'\circ\phi_i^{-1}=\P$ for every $i\in\enne$ and 
\begin{align*}
  (\varphi_{\lambda_{k_i}}', \varphi_{\lambda_{j_i}}')
  :=(\varphi_{\lambda_{k_i}}, \varphi_{\lambda_{j_i}})\circ\phi_i \to (\varphi'_1, \varphi'_2)
  \qquad&\text{in } (C^0([0,T]; H)\cap L^2(0,T; V_1))^2\,,\quad\P'\text{-a.s.}\,,
\end{align*}
for some measurable random variables
\[
  (\varphi'_1, \varphi'_2):(\Omega', \cF')\to (C^0([0,T]; H)\cap L^2(0,T; V_1))^2\,.
\]
Similarly, we have
\begin{align*}
(\bu_{\lambda_{k_i}}',\bu_{\lambda_{j_i}}') 
&:=(\bu_{\lambda_{k_i}}, \bu_{\lambda_{j_i}})\circ\phi_i
\to (\bu_1', \bu_2') &&\text{in } L^p(0,T; U)^2\,,
\quad\P'\text{-a.s.}\,,\\
(I_{\lambda_{k_i}}', I_{\lambda_{j_i}}')&:=
  (I_{\lambda_{k_i}}, I_{\lambda_{j_i}})\circ\phi_i\to (I_1', I_2')
   &&\text{in } C^0([0,T]; H)^2\,,
  \quad\P'\text{-a.s.}\,,\\
  W_i'&:=W\circ\phi_i \to W' &&\text{in } C^0([0,T]; K)\,, \quad\P'\text{-a.s.}
\end{align*}
for some measurable random variables 
\[
  (\bu_1', \bu_2'):(\Omega', \cF')\to L^p(0,T; U)
\]
and
\[
  (I_1', I_2'):(\Omega', \cF')\to C^0([0,T]; H)^2\,, \qquad
  W':(\Omega', \cF')\to C^0([0,T]; U)\,.
\]
Now, since $\bu_\lambda\to\bu$ in $L^p(0,T; U)$ $\P$-almost surely
on the whole sequence $\lambda$, for every 
arbitrary $f\in C^0(\erre)\cap L^\infty(\erre)$ we have
\begin{align*}
  {\E}'\left[f\left(\norm{\bu_1'-\bu_2'}_{L^p(0,T; U)}\right)\right]&=
  \lim_{i\to\infty}{\E}'\left[
  f\left(\norm{\bu_{\lambda_{k_i}}'-\bu_{\lambda_{j_i}}'}_{L^p(0,T; U)}\right)\right]\\
  &=
  \lim_{i\to\infty}\E\left[
  f\left(\norm{\bu_{\lambda_{k_i}}-\bu_{\lambda_{j_i}}}_{L^p(0,T; U)}\right)\right]=0\,,
\end{align*}
from which $\bu_1'=\bu_2'$ $\P'$-almost surely due to the arbitrariness of $f$.
Let us set then $\bu':=\bu_1'=\bu_2'$
and
$(\mu_{\lambda_{k_i}}',\mu_{\lambda_{j_i}}') 
:=(\mu_{\lambda_{k_i}}, \mu_{\lambda_{j_i}})\circ\phi_i$: since the 
maps $\phi_i$ preserve the laws, from the uniform estimates \eqref{est1}--\eqref{est3}
we deduce also that 
\begin{align*}
  (\varphi_{\lambda_{k_i}}',\varphi_{\lambda_{j_i}}')  \to (\varphi_1' , \varphi_2')
  \qquad&\text{in } L^\ell_\cP(\Omega'; C^0([0,T]; H)\cap L^2(0,T; V_1))^2
  \quad\forall\,\ell\in[1,p)\,,\\
  (\varphi_{\lambda_{k_i}}',\varphi_{\lambda_{j_i}}')  \wstarto (\varphi_1' , \varphi_2')
  \qquad&\text{in } L^p_w(\Omega'; L^\infty(0,T; V_1))^2\cap
  L^p_\cP(\Omega'; L^2(0,T; V_2))^2\,,\\
  (\mu_{\lambda_{k_i}}', \mu_{\lambda_{j_i}}' ) \wto (\mu_1', \mu_2') 
  \qquad&\text{in } L^{p/2}_\cP(\Omega'; L^2(0,T; V_1))^2\,,\\
  (\bu_{\lambda_{k_i}}', \bu_{\lambda_{j_i}}' ) \wstarto (\bu', \bu') 
  \qquad&\text{in } L^{\infty}_\cP(\Omega'; L^p(0,T; U))^2\,,
\end{align*}
for some measurable random variables
\[
  (\mu'_1, \mu'_2):(\Omega', \cF')\to  L^2(0,T; V_1)^2\,.
\]
Now, if we introduce the filtration
$(\cF_{i,t}')_{t\in[0,T]}$ as
\[
  \cF_{i,t}':=\sigma\{\varphi_{\lambda_{k_i}}'(s), \varphi_{\lambda_{j_i}}'(s),
  \mu_{\lambda_{k_i}}'(s),\mu_{\lambda_{j_i}}'(s),
  \bu_{\lambda_{k_i}}'(s), \bu_{\lambda_{j_i}}'(s),
  W_i'(s), I_{\lambda_{k_i}}'(s), I_{\lambda_{j_i}}'(s): s\leq t\}\,,
\]
using classical representation theorems for martingales
(see \cite{fland-gat} and \cite[\S~8.4]{dapratozab}) we have that 
$W_i'$ is a cylindrical Wiener process on $(\Omega', \cF', (\cF_t')_{t\in[0,T]}, \P')$ and
\[
  I_{\lambda_{k_i}}'=\int_0^\cdot 
  B(\varphi_{\lambda_{k_i}}'(s))\,\d W_i'(s)\,, \qquad
  I_{\lambda_{j_i}}'=\int_0^\cdot 
  B(\varphi_{\lambda_{j_i}}'(s))\,\d W_i'(s)\,,
\]
so that on the new probability space $(\Omega', \cF', \P')$ we have 
\begin{align*}
  \d\varphi_{\lambda_{k_i}}'
  -\Delta\mu_{\lambda_{k_i}}' \,\d t
  +\bu'_{\lambda_{k_i}}\cdot\nabla\varphi_{\lambda_{k_i}}' \,\d t =
  B(\varphi_{\lambda_{k_i}}') \,\d W_i'\,, \qquad&\varphi_{\lambda_{k_i}}'(0)=\varphi_0\,,\\
   \d\varphi_{\lambda_{j_i}}'
  -\Delta\mu_{\lambda_{j_i}}' \,\d t
  +\bu'_{\lambda_{j_i}}\cdot\nabla\varphi_{\lambda_{j_i}}'\,\d t =
  B(\varphi_{\lambda_{j_i}}' )\,\d W_i'\,,
  \qquad&\varphi_{\lambda_{j_i}}'(0)=\varphi_0\,,
\end{align*}
where the equations are intended in the usual variational sense \eqref{eq_ab_lam}.
Now, the strong convergences of $(\varphi_{\lambda_{k_i}}', \varphi_{\lambda_{j_i}}')_i$
imply, together with the Lipschitz-continuity of $B$, that
\[
  (B(\varphi_{\lambda_{k_i}}'), B(\varphi_{\lambda_{j_i}}'))
  \to (B(\varphi_1'), B(\varphi_2')) \qquad\text{in } 
  L^\ell_\cP(\Omega'; C^0([0,T]; \cL^2(K,H)))^2
  \quad\forall\,\ell\in[1,p)\,.
\]
Introducing then the limiting filtration $(\cF_{t}')_{t\in[0,T]}$ as
\[
  \cF_{t}':=\sigma\{\varphi_1'(s), \varphi_2'(s),
  \mu_1'(s), \mu_2'(s), \bu'(s),
  W'(s), I_1'(s), I_2'(s): s\leq t\}\,, \qquad t\in[0,T]\,,
\]
a classical argument based again on the martingale representation 
theorem (see \cite{fland-gat} and \cite[\S~8.4]{dapratozab}) yields the identification
\[
  I_1'=\int_0^\cdot B(\varphi_1'(s))\,\d W'(s)\,, \qquad
  I_2'=\int_0^\cdot B(\varphi_2'(s))\,\d W'(s)\,.
\]
Moreover, the strong convergences of $(\varphi_{\lambda_{k_i}}', \varphi_{\lambda_{j_i}}')_i$
together with the uniform estimate \eqref{est3} on the nonlinearities 
also give
\[
  (\Psi'_{\lambda_{k_i}}(\varphi_{\lambda_{k_i}}'), 
  \Psi'_{\lambda_{j_i}}(\varphi_{\lambda_{j_i}}'))\wto 
  (\Psi'(\varphi_1'), \Psi'(\varphi_2')) \qquad\text{in }
  L^{p/2}_\cP(\Omega'; L^2(0,T; H))^2\,.
\]
Putting all this information together, we deduce that 
$(\varphi_1', \varphi_2')$ solves the limit problem \eqref{eq1}--\eqref{eq4}
in the sense of Theorem~\ref{thm1} on the new probability space $(\Omega', \cF', \P')$, namely
\begin{align*}
  \d\varphi_1'
  -\Delta\mu_1' \,\d t
  +\bu'\cdot\nabla\varphi_1' \,\d t =
  B(\varphi_1') \,\d W'\,, \qquad&\varphi_1'(0)=\varphi_0\,,\\
   \d\varphi_2'
  -\Delta\mu_2' \,\d t
  +\bu'\cdot\nabla\varphi_2'\,\d t =
  B(\varphi_2' )\,\d W'\,,
  \qquad&\varphi_2'(0)=\varphi_0\,.
\end{align*}
Since we have already proved uniqueness of solutions in Subsection~\ref{ssec:uniq},
we deduce that 
\[
  \nu(\{(z_1,z_2)\in\mathcal X^2: z_1=z_2\})
  =\P'\left\{\norm{\varphi_1'-\varphi_2'}_{C^0([0,T]; H)\cap L^2(0,T; V_1)}=0\right\}=1\,.
\]
so that Lemma~\ref{lem:gyo-kry} ensures the strong convergence \eqref{strong}
also on the original probability space $(\Omega, \cF, \P)$.
Proceeding now in exactly the same way on $(\Omega, \cF, \P)$ instead, 
it is a standard matter to show that $(\varphi, \mu)$ is the unique 
solution to the state system \eqref{eq1}--\eqref{eq4}. 
Clearly, the global estimate \eqref{bound} follows directly by the computations in 
Subsection~\ref{ssec:est} and assumption {\bf A3}, 

\subsection{Continuous dependence}
Here we conclude the proof of Theorem~\ref{thm1} by showing the 
continuous dependence estimates \eqref{cont}--\eqref{cont2}.

First of all, \eqref{cont} is a consequence of the already proved
\eqref{bound} and Subsection~\ref{ssec:uniq}. 
Now, let us focus on proving \eqref{cont2}. To this end, 
we use the same notation of Subsection~\ref{ssec:uniq} and use It\^o's formula 
for the square of the $H$-norm instead, getting
\begin{align*}
  &\frac12\norm{\varphi(t)}_H^2 +
  \int_{Q_t}|\Delta\varphi|^2 - \int_{Q_t}(\Psi'(\varphi_1)-\Psi'(\varphi_2))\Delta\varphi
  +\int_{Q_t}\left(\bu\cdot\nabla\varphi_1 + \bu_2\cdot\nabla\varphi \right)\varphi\\
  &=\frac12\int_0^t\norm{B(\varphi_1(s))-B(\varphi_2(s))}^2_{\cL^2(K,H)}\,\d s
  +\int_0^t\left(\varphi(s),(B(\varphi_1(s))-B(\varphi_2(s)))\,\d W(s)\right)_H\,.
\end{align*}
The third term on the left hand side can be handled 
thanks to assumption {\bf A1}, the H\"older and Young inequalities,
and the embedding $V_1\embed L^6(\OO)$, as
\begin{align*}
  \int_{Q_t}(\Psi'(\varphi_1)-\Psi'(\varphi_2))\Delta\varphi
  &\leq\frac12\int_{Q_t}|\Delta\varphi|^2+
  c\int_{Q_t}\left(1+|\varphi_1|^4 + |\varphi_2|^4\right)|\varphi|^2\\
  &\leq\frac12\int_{Q_t}|\Delta\varphi|^2+
  \int_0^t\left(1+\norm{\varphi_1(s)}^4_{L^6(D)} + \norm{\varphi_2(s)}^4_{L^6(D)}\right)
  \norm{\varphi(s)}^2_{L^6(D)}\,\d s\\
  &\leq\frac12\int_{Q_t}|\Delta\varphi|^2+
  \left(1+\norm{\varphi_1}^4_{L^\infty(0,T; V_1)} + \norm{\varphi_2}^4_{L^\infty(0,T; V_1)}\right)
  \norm{\varphi}_{L^2(0,T; V_1)}^2\,.
  \end{align*}
The convection terms on the right-hand side can be treated similarly using the 
divergence theorem, 
the H\"older and Young inequalities, and the inclusion $L^6(\Omega)\embed V_1$ as
\[
  \int_{Q_t}\left(\bu\cdot\nabla\varphi_1 + \bu_2\cdot\nabla\varphi \right)\varphi=
  \int_{Q_t}\bu\cdot\nabla\varphi_1\varphi
  \leq
  \norm{\varphi}_{L^2(0,T; V_1)}^2+
  c\norm{\varphi_1}_{L^\infty(0,T; V_1)}^2\norm{\bu}_{L^2(0,T; U)}^2\,.
\]
Hence, we rearrange the terms and take power $p/6$ at both sides, obtaining,
thanks to the H\"older  and  Young inequalities, 
\begin{align*}
  &\E\norm{\varphi}_{L^\infty(0,T;H)}^{p/3} +
  \E\norm{\Delta\varphi}_{L^2(0,T; H)}^{p/3}\\
  &\leq c\left[1 + \norm{\varphi_1}_{L^p(\Omega; L^\infty(0,T; V_1))}^{\frac23p}
   +\norm{\varphi_2}_{L^p(\Omega; L^\infty(0,T; V_1))}^{\frac23p}\right]
   \norm{\varphi}_{L^p_\cP(\Omega;L^2(0,T; V_1))}^{p/3}
   +c\E\norm{\varphi}_{L^2(0,T; V_1)}^{p/3}\\
  &\qquad
  +c\norm{\varphi_1}_{L^p(\Omega; L^\infty(0,T; V_1))}^{p/3}\norm{\bu}_{\mU}^{p/3}
  +c\E\sup_{t\in[0,T]}\left|
  \int_0^t\left(\varphi(s),(B(\varphi_1(s))-B(\varphi_2(s)))\,\d W(s)\right)_H\right|^{p/6}
\end{align*}
where the Burkholder-Davis-Gundy inequality and the Lipschitz-continuity of $B$ yield
\[
\E\sup_{t\in[0,T]}\left|
  \int_0^t\left(\varphi(s),(B(\varphi_1(s))-B(\varphi_2(s)))\,\d W(s)\right)_H\right|^{p/6}
  \leq \sigma\E\norm{\varphi}_{L^\infty(0,T; H)}^{p/3}
  +c_\sigma\E\norm{\varphi}_{L^2(0,T; H)}^{p/3}
\]
for all $\sigma>0$. Hence, choosing $\sigma$ sufficiently small
and rearranging the terms, 
the continuous dependence \eqref{cont2} follows from
the already proved estimates \eqref{bound}--\eqref{cont}.
This concludes the proof of Theorem~\ref{thm1}.

%%%%%%%%%%%%%%%%%%%%%%%%%%%%%%%%%%%%

\section{Existence of optimal controls}
\label{sec:ex_optt}

In this section we prove Theorem~\ref{thm2}, showing that 
the optimisation problem {\bf(CP)} always admits a 
relaxed optimal control
 $\bu\in\mU_{ad}$ and a deterministic optimal control $\bu^{det}\in\mU_{ad}^{det}$.
The main idea idea is to use the direct method from calculus of variations, 
combined with a stochastic compactness argument.

Let $(\bu_n)_n\subset\mU_{ad}$ be a minimising sequence 
for the functional $\widetilde J$, in the sense that 
\[
  \tilde J(\bu_n)\searrow \inf_{\bv\in\mU_{ad}}\widetilde J(\bv)\,,
\]
and define $(\varphi_n,\mu_n)_n$ as the unique respective 
solutions to the state system \eqref{eq1}--\eqref{eq4},
in the sense of Theorem~\ref{thm1}.
Thanks to the definition of $\mU_{ad}$
and the estimate \eqref{bound}, we deduce that 
there exist $\bu\in\mU_{ad}$ and a triplet $(\varphi,\mu,\xi)$ with 
\begin{align*}
  &\varphi \in L^p_\cP(\Omega; C^0([0,T]; H)\cap L^2(0,T; V_2))\cap 
  L^p_w(\Omega; L^\infty(0,T; V_1))\,,\\
  &\mu\in L^{p/2}_\cP(\Omega; L^2(0,T; V_1))\,, \qquad
  \xi \in L^{p/2}_\cP(\Omega; L^2(0,T; H))\,,
\end{align*}
such that, as $n\to\infty$, possibly on a subsequence, 
\begin{align*}
  \varphi_n\wstarto \varphi\qquad&\text{in }
  L^p_w(\Omega; L^\infty(0,T; V_1))\cap L^p_\cP(\Omega; L^2(0,T; V_2))\,,\\
  \mu_n\wto\mu\qquad&\text{in }
  L^{p/2}_\cP(\Omega; L^2(0,T; V_1))\,,\\
  \Psi'(\varphi_n)\wto\xi\qquad&\text{in }
  L^{p/2}_\cP(\Omega; L^2(0,T; H))\,,\\
  \bu_n\wstarto\bu\qquad&\text{in }
  L^{\infty}_\cP(\Omega; L^p(0,T; U))\,.
\end{align*}
Assumption {\bf A3} and the uniform estimates on $(\varphi_n)_n$ ensure also that 
\[
  \norm{B(\varphi_n)}_{L^\infty(\Omega\times(0,T); \cL^2(K,H))
  \cap L^p(\Omega; L^\infty(0,T; \cL^2(K,V_1)))} \leq  c\,,
\]
so that in particular 
\[
  \norm{I_n:=\int_0^\cdot B(\varphi_n(s))\,\d W(s)}_{L^p_\cP(\Omega; W^{s,p}(0,T; V_1))} \leq c_s
  \quad\forall\,s\in(0,1/2)\,.
\]
By comparison in the equation \eqref{eq1} we infer then
\[
  \norm{\varphi_n}_{L^p_\cP(\Omega; W^{s,p}(0,T; V_1^*))} \leq c_s
  \quad\forall\,s\in(0,1/2)\,,
\]
which ensures that the laws of $(\varphi_n)_n$ are tight
on the space $C^0([0,T]; H)\cap L^2(0,T; V_1)$.
We argue now on the same line of Subsection~\ref{ssec:pass_lim}.
As a consequence of the Skorokhod theorem there is a probability space 
$(\Omega', \cF', \P')$ and measurable maps 
$\phi_i:(\Omega', \cF')\to (\Omega, \cF)$
with $\P'\circ\phi_i^{-1}=\P$ for all $i\in\enne$, such that
\begin{align*}
  \varphi_{n_i}':=\varphi_{n_i}\circ\phi_i \to \varphi'
  \qquad&\text{in } L^\ell_\cP(\Omega'; C^0([0,T]; H)\cap L^2(0,T; V_1))
  \quad\forall\,\ell\in[1,p)\,,\\
  \varphi_{n_i}' \wstarto \varphi'
  \qquad&\text{in } L^p_w(\Omega'; L^\infty(0,T; V_1))\cap
  L^p_\cP(\Omega'; L^2(0,T; V_2))\,,\\
  \mu_{n_i}' :=\mu_{n_i}\circ\phi_i\wto \mu'
  \qquad&\text{in } L^{p/2}_\cP(\Omega'; L^2(0,T; V_1))\,,\\
  \bu_{n_i}':=\bu_{n_i}\circ\phi_i\wstarto\bu'
  \qquad&\text{in } L^{\infty}_\cP(\Omega'; L^p(0,T; U))\,,\\
  \varphi'_{Q,i}:=\varphi_Q\circ\phi_i'\wto\varphi_Q'
  \qquad&\text{in } L^2_\cP(\Omega'; L^2(0,T; H))\,,\\
  \varphi_{T,i}':=\varphi_T\circ\phi_i \wto \varphi_T'
  \qquad&\text{in } L^2(\Omega',\cF_T'; H)\,.
\end{align*}
Furthermore, on the new probability space we have
\begin{align*}
  \d\varphi_{n_i}'
  -\Delta\mu_{n_i}' \,\d t
  +\bu_{n_i}'\cdot\nabla\varphi_{n_i}' \,\d t =
  B(\varphi_{n_i}') \,\d W_i'\,, \qquad&\varphi_{n_i}'(0)=\varphi_0\,,
\end{align*}
where the stochastic integral are intended with respect to 
a suitably defined filtration $(\cF_{i,t})_{t\in[0,T]}$.
Proceeding as in Subsection~\ref{ssec:pass_lim}, we infer that 
\[
  \Psi'(\varphi_{n_i}')\wto 
  \Psi'(\varphi') \qquad\text{in }
  L^{p/2}_\cP(\Omega'; L^2(0,T; H))\,,
\]
so that by assumption {\bf A3} and the 
martingale representation theorem we 
can pass to the limit as $i\to \infty$ on the new 
probability space and get 
\begin{align*}
  \d\varphi'
  -\Delta\mu' \,\d t
  +\bu'\cdot\nabla\varphi' \,\d t =
  B(\varphi') \,\d W'\,, \qquad&\varphi'(0)=\varphi_0\,.
\end{align*}
This shows that $\bu'\in\mU_{ad}'$ and that $(\varphi',\mu')=S'(\bu')$.
To conclude that $\bu'$ is a relaxed 
optimal control for the optimisation problem {\bf (CP)},
we note that by the weak lower semicontinuity of
the cost functional $J$ we have
\begin{align*}
  \widetilde J'(\bu')&=
  \frac{\alpha_1}{2}{\E}'\int_Q|\varphi'-\varphi'_Q|^2 + 
  \frac{\alpha_2}{2}{\E}'\int_\OO|\varphi'(T)-\varphi'_T|^2 +
  \frac{\alpha_3}{2}\E\int_Q|\bu'|^2\\
  &\leq\liminf_{i\to\infty}
  \left(
  \frac{\alpha_1}{2}{\E}'\int_Q|\varphi_{n_i}'-\varphi_{Q,i}'|^2 + 
  \frac{\alpha_2}{2}{\E}'\int_\OO|\varphi_{n_i}'(T)-\varphi_{T,i}'|^2 +
  \frac{\alpha_3}{2}{\E}'\int_Q|\bu_{n_i}'|^2\right)\\
  &=\liminf_{i\to\infty}
  \left(
  \frac{\alpha_1}{2}\E\int_Q|\varphi_{n_i}-\varphi_Q|^2 + 
  \frac{\alpha_2}{2}\E\int_\OO|\varphi_{n_i}(T)-\varphi_T|^2 +
  \frac{\alpha_3}{2}\E\int_Q|\bu_{n_i}|^2\right)\\
  &=\liminf_{n\to\infty}\widetilde J(\bu_n)=
  \inf_{\bv\in\mathcal U_{ad}}\tilde J(\bv)\,,
\end{align*}
so that $\bu'\in\mU'_{ad}$ is a relaxed optimal control in the sense of Definition~\ref{def:OC}.

In order to show existence of a deterministic optimal control, 
the argument is similar. We start taking a minimising sequence 
$(\bu_n)_n\subset\mU_{ad}^{det}$ such that
\[
  \tilde J(\bu_n)\searrow \inf_{\bv\in\mU^{det}_{ad}}\widetilde J(\bv)\,.
\]
Arguing exactly as above, thanks to the fact that $(\bu_n)_n$ are
deterministic, in this case we have that 
$\bu_{n_i}'=\bu_{n_i}$ for every $i\in\enne$. Consequently, 
in this case we can
$(\varphi_n)_n$ inherits some strong compactness 
properties on the original probability space,
using a similar argument to the one of Subsection~\ref{ssec:pass_lim},
by employing Lemma~\ref{strong}.
Namely, we infer the strong convergence
\[
  \varphi_n\to \varphi \qquad\text{in } C^0([0,T];H)\cap L^2(0,T; V_1)\,, \quad\P\text{-a.s.}
\]
on the original probability space $(\Omega,\cF,\P)$. It follows then that 
$\xi=\Psi'(\varphi)$ almost everywhere, and letting $n\to\infty$ yields 
\[
   \d\varphi
  -\Delta\mu\,\d t
  +\bu\cdot\nabla\varphi \,\d t =
  B(\varphi) \,\d W\,, \qquad\varphi(0)=\varphi_0\,,
\]
so that $(\varphi,\mu)=S(\bu)$. At this point, the conclusion
follows as above by lower semicontinuity of the cost functional.

%%%%%%%%%%%%%%%%%%%%%%%%%%%%%%%%%

\section{Linearised system and differentiability of the control-to-state map}
\label{sec:lin}
The aim of this section is to prove that the linearised state system \eqref{eq1_lin}--\eqref{eq2_lin}
is well-posed and to characterise its solution as the derivative on the 
control-to-state map. Namely, we prove here Theorem~\ref{thm3}.

\subsection{Existence}
\label{ssec:ex_lin}
Let $\bu\in\widetilde \mU_{ad}$ and $\bh\in\mU$ be arbitrary and fixed.
Using the notation of Subsection~\ref{ssec:approx}, we consider the approximated 
linearised problem
\begin{align}
    \label{eq1_lin_lam'}
    \d\theta_{\bh,\lambda} - \Delta \nu_{\bh,\lambda}\,\d t 
    +\bh\cdot\nabla\varphi\,\d t
    + \bu_\lambda\cdot\nabla\theta_{\bh,\lambda}\,\d t
    = DB(\varphi)\theta_{\bh,\lambda}\,\d W \qquad&\text{in } (0,T)\times \OO\,,\\
    \label{eq2_lin_lam'}
    \nu_{\bh,\lambda}=-\Delta \theta_{\bh,\lambda} +
     \Psi_\lambda''(\varphi)\theta_{\bh,\lambda}
     \qquad&\text{in } (0,T)\times \OO\,,\\
    \label{eq3_lin_lam'}
    {\bf n}\cdot\nabla\theta_{\bh,\lambda} = {\bf n}\cdot\nabla\nu_{\bh,\lambda} = 0 
    \qquad&\text{in } (0,T)\times\partial \OO\,,\\
    \label{eq4_lin_lam'}
     \theta_{\bh,\lambda}(0)=0 \qquad&\text{in } \OO\,.
\end{align}
Noting that $\Psi''_\lambda(\varphi)\in L^\infty(\Omega\times Q)$, 
the classical variational approach ensures existence and uniqueness of 
the approximated solution
\begin{align*}
  &\theta_{\bh,\lambda} \in 
  L^2_{\cP}\left(\Omega; C^0([0,T]; H)\cap L^2(0,T; V_2)\right)\,,\\
  &\nu_{\bh,\lambda}=-\Delta\theta_{\bh,\lambda}
  +\Psi''_\lambda(\varphi)\theta_{\bh,\lambda}\in L^2_{\cP}(\Omega; L^2(0,T; H))\,,
\end{align*}
in the sense that, for every $\zeta\in V_2$, for every $t\in[0,T]$, $\P$-almost surely,
\begin{align}
   \label{var_lin_app'}
    &\left(\theta_{\bh,\lambda}(t),\zeta\right)_H 
    - \int_{Q_t}\nu_{\bh,\lambda}\Delta\zeta
  -\int_{Q_t}(\varphi\bh + \theta_{\bh,\lambda}\bu)\cdot\nabla\zeta
  = \left(\int_0^tDB(\varphi)\theta_{\bh,\lambda}\,\d W(s), \zeta\right)_H\,.
\end{align}

Noting that $(\theta_{\bh,\lambda})_\OO=0$, we can write It\^o's formula 
for $\frac12\norm{\nabla\mN\theta_{\bh,\lambda}}_H^2$,
getting
\begin{align*}
  &\frac12\norm{\nabla\mN\theta_{\bh,\lambda}(t)}_H^2
  +\int_{Q_t}|\nabla\theta_{\bh,\lambda}|^2
  =-\int_{Q_t}\Psi_\lambda''(\varphi)|\theta_{\bh,\lambda}|^2
  +\int_{Q_t}\left(\varphi\bh+\theta_{\bh,\lambda}\bu_\lambda\right)
  \cdot\nabla\mN\theta_{\bh,\lambda}\\
  &\qquad+\frac12\int_0^t\norm{\nabla\mN DB(\varphi(s))\theta_\bh(s)}_{\cL^2(K,H)}^2\,\d s
  +\int_0^t\left(\mN\theta_{\bh,\lambda}(s), 
  DB(\varphi(s))\theta_{\bh,\lambda}(s)\,\d W(s)\right)_H\,.
\end{align*}
Now, assumption {\bf A1}, the H\"older--Young inequalities and the compactness inequality 
\eqref{comp_ineq},
and the embedding $V_1\embed L^6(\OO)$ give, for all $\eps>0$,
\begin{align*}
  &-\int_{Q_t}\Psi_\lambda''(\varphi)|\theta_{\bh,\lambda}|^2
  +\int_{Q_t}\left(\varphi\bh+\theta_{\bh,\lambda}\bu_\lambda\right)
  \cdot\nabla\mN\theta_{\bh,\lambda}\\
  &\leq\eps\int_{Q_t}|\nabla\theta_{\bh,\lambda}|^2+
  \norm{\varphi}_{L^\infty(0,T; V_1)}^2+
  c_\eps\int_0^t\left(1+\norm{\bh(s)}_U^2+\norm{\bu_\lambda(s)}_U^2\right)
  \norm{\nabla\mN\theta_{\bh,\lambda}(s)}_H^2\,\d s\,.
\end{align*}
Similarly, by {\bf C2} and again the compactness inequality 
\eqref{comp_ineq}, we have 
\[
  \frac12\int_0^t\norm{\nabla\mN DB(\varphi(s))\theta_\bh(s)}_{\cL^2(K,H)}^2\,\d s\leq
  \eps\int_{Q_t}|\nabla\theta_{\bh,\lambda}|^2 + 
  c_\eps\int_0^t
  \norm{\nabla\mN\theta_{\bh,\lambda}(s)}_H^2\,\d s\,.
\]
As for the stochastic integral, 
the Burkholder-Davis-Gundy and Young inequalities give
(see for example \cite[Lem.~4.1]{mar-scar-ref}),
together with \eqref{comp_ineq} and {\bf C2}
\begin{align*}
  &\E\sup_{r\in[0,t]}\left|\int_0^r\left(\mN\theta_{\bh,\lambda}(s), 
  DB(\varphi(s))\theta_{\bh,\lambda}(s)\,\d W(s)\right)_H\right|^{p/2}\\
  &\leq \eps\E\norm{\mN\theta_{\bh,\lambda}}_{L^\infty(0,t; H)}^p
  +c_\eps\E\norm{\theta_{\bh,\lambda}}_{L^2(0,t; H)}^p\\
  &\leq\eps\E\norm{\nabla\mN\theta_{\bh,\lambda}}_{L^\infty(0,t; H)}^p+
  \eps\E\norm{\nabla\theta_{\bh,\lambda}}_{L^\infty(0,t; H)}^p
  +c_\eps\E\norm{\nabla\mN\theta_{\bh,\lambda}}_{L^2(0,t; H)}^p\,.
\end{align*}
Consequently, using the same iterative--patching argument of
Subsection~\ref{ssec:uniq},
raising to power $p/2$, taking supremum in time and expectations, 
we infer that 
\beq
  \label{est1_lin}
  \norm{\theta_{\bh,\lambda}}_{L^p_\cP(\Omega; C^0([0,T]; V_1^*)\cap L^2(0,T; V_1))} \leq c\,.
\eeq
Now, It\^o's formula for $\frac12\norm{\theta_{\bh,\lambda}}_H^2$ yields 
\begin{align*}
  &\frac12\norm{\theta_{\bh,\lambda}}_H^2
  +\int_{Q_t}|\Delta\theta_{\bh,\lambda}|^2
  =\int_{Q_t}\left(\varphi\bh+\theta_{\bh,\lambda}\bu_\lambda\right)
  \cdot\nabla\theta_{\bh,\lambda}
  +\int_{Q_t}\Psi''_\lambda(\varphi)\theta_{\bh,\lambda}\Delta\theta_{\bh,\lambda}\\
  &\qquad+\frac12\int_0^t
  \norm{DB(\varphi(s))\theta_{\bh,\lambda}(s)}_{\cL^2(K,H)}^2\,\d s
  +\int_0^t\left(\theta_{\bh,\lambda}(s), 
  DB(\varphi(s))\theta_{\bh,\lambda}(s)\,\d W(s)\right)_H\,,
\end{align*}
where by the divergence theorem we have 
\[
  \int_{Q_t}\left(\varphi\bh+\theta_{\bh,\lambda}\bu_\lambda\right)
  \cdot\nabla\theta_{\bh,\lambda} =
  \int_{Q_t}\varphi\bh
  \cdot\nabla\theta_{\bh,\lambda}\,.
\]
Hence, it is not difficult to see that,
using again
the H\"older, Young and Burkholder--Davis--Gundy inequalities, 
assumption {\bf C2}, and the estimate \eqref{est1_lin},
all the terms on the right-hand side can be handled, 
except the one containing $\Psi''$. For this one, we proceed 
using {\bf C1}, the embedding $V_1\embed L^6(\OO)$, as
\begin{align*}
  \int_{Q_t}\Psi''_\lambda(\varphi)\theta_{\bh,\lambda}\Delta\theta_{\bh,\lambda}
  &\leq\eps\int_{Q_t}|\Delta\theta_{\bh,\lambda}|^2
  +c_\eps\int_0^1\left(1 + \norm{\varphi(s)}_{V_1}^4\right)
  \norm{\theta_{\bh,\lambda}(s)}_{V_1}^2\,\d s\\
  &\leq\eps\int_{Q_t}|\Delta\theta_{\bh,\lambda}|^2
  +c_\eps\left(1 + \norm{\varphi}_{L^\infty(0,T; V_1)}^4\right)
  \norm{\theta_{\bh,\lambda}}_{L^2(0,T; V_1)}^2\,,
\end{align*}
where, thanks to \eqref{est1_lin} and the H\"older inequality, 
\[
  \norm{\norm{\varphi}_{L^\infty(0,T; V_1)}^4
  \norm{\theta_{\bh,\lambda}}_{L^2(0,T; V_1)}^2}_{L^{p/6}(\Omega)}
  \leq\norm{\varphi}_{L^p(\Omega; L^\infty(0,T; V_1))}^4
  \norm{\theta_{\bh,\lambda}}_{L^p_\cP(\Omega; L^2(0,T; V_1))}^2\leq c\,.
\]
Consequently, we deduce that
\beq
  \label{est2_lin}
  \norm{\theta_{\bh,\lambda}}_{L^{p/3}_\cP(\Omega; C^0([0,T]; H)\cap L^2(0,T; V_2))} \leq c\,,
\eeq
from which, by comparison in \eqref{eq2_lin_lam'},
\beq
  \label{est3_lin}
  \norm{\nu_{\bh,\lambda}}_{L^{p/3}_\cP(\Omega;  L^2(0,T; H))} \leq c\,,
\eeq
We infer the existence of $(\theta_\bh, \nu_\bh)$ with
\begin{align*}
  &\theta_\bh \in 
  L^{p}_{\cP}\left(\Omega; C^0([0,T]; V_1^*)\cap L^2(0,T; V_1)\right)
  \cap L^{p/3}_w(\Omega; L^\infty(0,T; H))\cap L^{p/3}_\cP(\Omega; L^2(0,T; V_2))\,,\\
  &\nu_\bh \in L^{p/3}_{\cP}(\Omega; L^2(0,T; H))\,,
\end{align*}
such that, as $\lambda\searrow0$ (possibly on a subsequence),
\begin{align}
  \label{conv1_lin'}
  \theta_{\bh,\lambda}\wstarto \theta_\bh
  \qquad&\text{in }  L^{p}_{w}\left(\Omega; L^\infty(0,T; V_1^*)\right)\cap
  L^{p/3}_{w}\left(\Omega; L^\infty(0,T; H)\right)\,,\\
  \label{conv2_lin'}
  \theta_{\bh,\lambda}\wto \theta_\bh
  \qquad&\text{in }  L^{p}_\cP\left(\Omega; L^2(0,T; V_1)\right)\cap
  L^{p/3}_\cP\left(\Omega; L^2(0,T; V_2)\right)\,,\\
  \label{conv3_lin'}
  \nu_{\bh,\lambda}\wto \nu_\bh
  \qquad&\text{in } 
  L^{p/3}_\cP\left(\Omega; L^2(0,T; H)\right)\,.
\end{align}
Since the systems \eqref{eq1_lin_lam'}--\eqref{eq4_lin_lam'}
and \eqref{eq1_lin}--\eqref{eq4_lin} are linear, the passage to the 
limit is straightforward. Indeed, 
by assumption {\bf C2} and the dominated convergence theorem, it follows that 
\[
  DB(\varphi)\theta_{\bh,\lambda}\wto DB(\varphi)\theta_\bh
  \qquad\text{in }  L^{p}_\cP\left(\Omega; L^2(0,T; \cL^2(K,H))\right)\,.
\]
Moreover, thanks to {\bf C1} and the regularity of $\varphi$ we have 
$\Psi''(\varphi) \in L^3(\Omega; L^\infty(0,T; L^3(\OO)))$, so 
in particular
\[
  \Psi''_\lambda(\varphi) \to \Psi''(\varphi) \qquad\text{in } L^3(\Omega\times Q)\,,
\]
and also, thanks to \eqref{conv2_lin'},
\[
  \Psi''_\lambda(\varphi)\theta_{\bh,\lambda}
  \wto \Psi''(\varphi)\theta_\bh 
  \qquad\text{in } L^{6/5}_\cP(\Omega; L^{6/5}(0,T; L^{6/5}(\OO)))\,.
\]
We deduce that letting $\lambda\searrow0$ in \eqref{var_lin_app'}
we get that $(\theta_\bh, \nu_\bh)$ is a solution to \eqref{eq1_lin}--\eqref{eq4_lin}
in the sense of Theorem~\ref{thm3}. The strong continuity in $H$ of $\theta_\bh$
follows a posteriori with a classical method by It\^o's formula on the limit equation \eqref{eq1_lin}.

\subsection{Uniqueness}
We show here that the linearised system
\eqref{eq1_lin}--\eqref{eq4_lin} admits at most one solution.
By linearity, it enough to check that if $(\theta,\nu)$ is a solution 
to \eqref{eq1_lin}--\eqref{eq4_lin} in the sense of Theorem~\ref{thm3}
with $\bh=0$, then $\theta=\nu=0$.
To this end, we note that \eqref{eq1_lin} yields $\theta_\OO=0$,
so that It\^o's formula gives
\begin{align*}
  \frac12\norm{\nabla\mN\theta(t)}_H^2 + \int_{Q_t}|\nabla\theta|^2+\int_{Q_t}\Psi''(\varphi)|\theta|^2
  &=\int_{Q_t}\theta\bu\cdot\nabla\mN\theta 
  +\int_0^t\left(\mN\theta(s), DB(\varphi(s))\theta(s)\right)_H\\
  &+ \frac12\int_0^t\norm{\nabla\mN DB(\varphi(s))\theta(s)}_{\cL^2(K,H)}^2\,\d s\,.
\end{align*}
Now, we can argue on the same line of Subsection~\ref{ssec:ex_lin}
by using assumption {\bf A1} on $\Psi''$, {\bf C2} on $DB$, together 
with Burkholder-Davis-Gundy and Young inequalities to get 
\[
  \norm{\theta}_{L^p_\cP(\Omega; C^0([0,T]; V_1^*)\cap L^2(0,T; V_1))}\leq 0\,,
\]
from which $\theta=0$, and also $\nu=0$ by comparison in \eqref{eq2_lin}.
This show that the linearised system \eqref{eq1_lin}--\eqref{eq4_lin}
admits at most one solution.

\subsection{G\^ateaux-differentiability}
We prove here that $S_1$ is G\^ateaux-differentiable.
Let $\bu\in\widetilde \mU_{ad}$ and $\bh\in\mU$ be arbitrary and fixed:
since $\widetilde \mU_{ad}$ is open in $\mU$, there exists $\delta_0>0$
such that $\bu+\delta\bh\in\widetilde \mU_{ad}$ for all $\delta\in[-\delta_0, \delta_0]$.
For every such $\delta$, setting $(\varphi_\delta, \mu_\delta):=S(\bu+\delta\bh)$
and $(\varphi, \mu):=S(\bu)$, the difference of the respective equations
(for $\delta\neq0$) gives
\begin{align*}
  &\d\left(\frac{\varphi_\delta-\varphi}{\delta}\right) 
  - \Delta\left(\frac{\mu_\delta-\mu}{\delta}\right)\,\d t
  +\bu\cdot\nabla\left(\frac{\varphi_\delta-\varphi}{\delta}\right)\,\d t
  +\bh\cdot\nabla\varphi_\delta\,\d t=
  \frac{B(\varphi_\delta)-B(\varphi)}{\delta}\,\d W\,,\\
  &\frac{\mu_\delta-\mu}{\delta} = 
  -\Delta\left(\frac{\varphi_\delta-\varphi}{\delta} \right)
  +\frac{\Psi'(\varphi_\delta)-\Psi'(\varphi)}{\delta}\,,\\
  &\left(\frac{\varphi_\delta-\varphi}{\delta}\right)(0)=0\,,
\end{align*}
whose natural variational formulation reads
\begin{align}
  \nonumber
    &\left(\frac{\varphi_\delta-\varphi}{\delta}(t),\zeta\right)_H 
    - \int_{Q_t}\frac{\mu_\delta-\mu}{\delta}\Delta\zeta
  -\int_{Q_t}(\varphi_\delta\bh + \frac{\varphi_\delta-\varphi}{\delta}\bu)\cdot\nabla\zeta\\
  \label{var_lin_app}
  &= \left(\int_0^t\frac{B(\varphi_\delta(s))-B(\varphi(s))}{\delta}\,\d W(s), \zeta\right)_H
  \qquad\forall\,\zeta\in V_2\,,\qquad
  \forall\,t\in[0,T]\,,\quad\P\text{-a.s.}
\end{align}
Now, by the continuous dependence estimate \eqref{cont2}, 
we deduce that there exists a constant $c>0$ independent of $\delta$ such that 
\begin{align*}
  &\norm{\frac{\varphi_\delta-\varphi}{\delta}}_{L^p_\cP(\Omega; C^0([0,T]; V_1^*)
  \cap L^2(0,T; V_1))}
  \leq c\,,\\
  &\norm{\frac{\varphi_\delta-\varphi}{\delta}}_{L^{p/3}_\cP(\Omega; C^0([0,T]; H)
  \cap L^2(0,T, V_2))}+
  \norm{\frac{\mu_\delta-\mu}{\delta}}_{L^{p/3}_\cP(\Omega; L^2(0,T; H))}\leq c\,,
\end{align*}
so that there exist $(\theta_\bh, \nu_\bh)$ with
\begin{align*}
  &\theta_\bh \in L^{p}_w(\Omega; L^\infty(0,T; V_1^*))\cap L^{p}_\cP(\Omega; L^2(0,T; V_1))\cap
  L^{p/3}_w(\Omega; L^\infty(0,T; H))\cap L^{p/3}_\cP(\Omega; L^2(0,T; V_2))\,,\\
  &\nu_\bh \in L^{p/3}_\cP(\Omega; L^2(0,T; H))\,,
\end{align*}
such that, as $\delta\to0$ possibly on a subsequence,
\begin{align}
  \label{conv1_lin}
  \frac{\varphi_\delta-\varphi}{\delta}\wstarto\theta_\bh
  \qquad&\text{in } 
  L^{p}_w(\Omega; L^\infty(0,T; V_1^*))\cap L^{p/3}_w(\Omega; L^\infty(0,T; H))\,,\\
  \label{conv2_lin}
  \frac{\varphi_\delta-\varphi}{\delta}\wto\theta_\bh
  \qquad&\text{in } 
  L^{p}_\cP(\Omega; L^2(0,T; V_1))\cap L^{p/3}_\cP(\Omega; L^2(0,T; V_2))\,,\\
  \label{conv3_lin}
  \frac{\mu_\delta-\mu}{\delta}\wto\nu_\bh
  \qquad&\text{in } 
  L^{p/3}_\cP(\Omega; L^2(0,T; H))\,.
\end{align}
It follows in particular that 
\beq\label{conv4_lin}
  \varphi_\delta\to\varphi \quad\text{in } 
  L^p_\cP(\Omega; C^0([0,T]; V_1^*)\cap L^2(0,T; V_1)) \cap 
  L^{p/3}_\cP(\Omega; C^0([0,T]; H)\cap L^2(0,T, V_2))\,.
\eeq
Furthermore, since $\bu\in \mU$, by the inclusion $V_1\embed L^6(\OO)$,
the H\"older inequality,
and the convergence \eqref{conv2_lin} it holds that
\beq\label{conv5_lin}
  \left(\frac{\varphi_\delta-\varphi}{\delta}\right)\bu\wto \theta_\bh\bu
  \qquad\text{in } L^p_\cP(\Omega; L^{\frac{2p}{p+2}}(0,T;H^d))\,.
\eeq
As far as the nonlinear term is concerned, thanks to the mean-value theorem we have
\begin{align*}
  &\frac{\Psi'(\varphi_\delta) - \Psi'(\varphi)}\delta - \Psi''(\varphi)\theta_\bh\\
  &=\frac{\Psi'(\varphi_\delta) - \Psi'(\varphi) - \Psi''(\varphi)(\varphi_\delta-\varphi)}\delta
  +\Psi''(\varphi)\left(\frac{\varphi_\delta - \varphi}{\delta} - \theta_\bh\right)\\
  &=\frac{\varphi_\delta - \varphi}{\delta}
  \int_0^1\left(\Psi''(\varphi + s(\varphi_\delta-\varphi)) - \Psi''(\varphi)\right)\,\d s
  +\Psi''(\varphi)\left(\frac{\varphi_\delta - \varphi}{\delta} - \theta_\bh\right)\,.
\end{align*}
Now, by the strong convergence \eqref{conv4_lin} and the continuity of $\Psi''$ we have 
\[
  \Psi''(\varphi + s(\varphi_\delta-\varphi)) - \Psi''(\varphi) \to 0 \quad\forall\,s\in[0,1]\,,\quad
  \text{a.e.~in } \Omega\times(0,T)\times \OO\,,
\]
where, recalling that by {\bf C1} $\Psi''$ has quadratic growth, 
thanks to the embedding $V_1\embed L^6(\OO)$
the left-hand side is uniformly bounded in the space $L^{p/2}(\Omega; L^\infty(0,T; L^3(\OO)))$,
so that
\[
  \int_0^1\left(\Psi''(\varphi + s(\varphi_\delta-\varphi)) - \Psi''(\varphi)\right)\,\d s
  \to 0 \qquad\text{in } L^{\ell'}_\cP(\Omega; L^{\ell''}(0,T; L^3(\OO)))
\]
for every $\ell'\in[1,p/2)$ and $\ell''\in[1,+\infty)$.
Taking \eqref{conv2_lin} into account, 
we infer in particular that 
\[
  \frac{\varphi_\delta - \varphi}{\delta}
  \int_0^1\left(\Psi''(\varphi + s(\varphi_\delta-\varphi)) - \Psi''(\varphi)\right)\,\d s \wto 0
  \qquad\text{in } L^{\ell'}_\cP(\Omega; L^{\ell''}(0,T; H))
\]
for every $\ell'\in[1,p/3)$ and $\ell''\in[1,2)$.
Similarly, thanks to {\bf C1} and the 
regularity of $\varphi$ we have
$\Psi''(\varphi)\in L^{p/2}(\Omega; L^\infty(0,T; L^3(\OO)))$,
and the same argument as above yields 
\[
  \Psi''(\varphi)\left(\frac{\varphi_\delta - \varphi}{\delta} - \theta_\bh\right) \wto 0
  \qquad\text{in } L^{\ell'}_\cP(\Omega; L^{\ell''}(0,T; H))
\]
for every $\ell'\in[1,p/3)$ and $\ell''\in[1,2)$.
It follows that
\beq\label{conv6_lin}
  \frac{\Psi'(\varphi_\delta) - \Psi'(\varphi)}\delta \wto 
  \Psi''(\varphi)\theta_\bh \quad\text{in } L^{\ell'}_\cP(\Omega; L^{\ell''}(0,T; H))
  \qquad\forall\,\ell'\in[1,p/3)\,,\quad\forall\,\ell''\in[1,2)\,.
\eeq
Lastly, let us handle the stochastic integral.
By the Lipschitz-continuity of $B$ in {\bf A3} we have
\begin{align*}
  &\frac{B(\varphi_\delta)-B(\varphi)}{\delta} - DB(\varphi)\theta_\bh\\
  &=\frac{B(\varphi_\delta)-B(\varphi) - DB(\varphi)(\varphi_\delta-\varphi)}{\delta} 
  +DB(\varphi)\left(\frac{\varphi_\delta-\varphi}{\delta} - \theta_\bh\right)\\
  &=\int_0^1\left(DB(\varphi+s(\varphi_\delta-\varphi))-DB(\varphi)\right)
  \frac{\varphi_\delta-\varphi}\delta\,\d s
  +DB(\varphi)\left(\frac{\varphi_\delta-\varphi}{\delta} - \theta_\bh\right)\,.
\end{align*}
Now, the strong convergence \eqref{conv4_lin}, the continuity
and boundedness of $DB$ in {\bf C2} imply together with the dominated 
convergence theorem that 
\[
  \int_0^1\left(DB(\varphi+s(\varphi_\delta-\varphi))-DB(\varphi)\right)\,\d s \to 0
  \qquad\text{in } L^\ell(\Omega; L^\ell(0,T; \cL(V_1,\cL^2(K,H))))
\]
for every $\ell\in[1,+\infty)$. Since
$\frac{\varphi_\delta-\varphi}{\delta}$ is bounded in $L^{p/3}(\Omega; L^4(0,T; V_1))$
by interpolation of \eqref{conv1_lin}--\eqref{conv2_lin}, it follows that 
\[
  \int_0^1\left(DB(\varphi+s(\varphi_\delta-\varphi))-DB(\varphi)\right)
  \frac{\varphi_\delta-\varphi}\delta\,\d s \wto 0
  \qquad\text{in } L^\ell(\Omega; L^2(0,T; \cL^2(K,H)))
\]
for every $\ell\in[1,p/3)$. Similarly, by the boundedness of 
$DB$ in {\bf C2} and the convergence \eqref{conv2_lin} we have also
\[
  DB(\varphi)\left(\frac{\varphi_\delta-\varphi}{\delta} - \theta_\bh\right)\wto0
  \qquad\text{in } L^p_\cP(\Omega; L^2(0,T; \cL^2(K,H)))\,.
\]
Hence, we obtain that 
\beq\label{conv7_lin}
  \frac{B(\varphi_\delta)-B(\varphi)}{\delta} \wto DB(\varphi)\theta_\bh
  \qquad\text{in } L^\ell(\Omega; L^2(0,T; \cL^2(K,H)))\quad
  \forall\,\ell\in[1,p/3)\,.
\eeq
Finally, letting $\delta\to0$ in \eqref{var_lin_app}
using convergences \eqref{conv1_lin}--\eqref{conv7_lin},
we deduce that actually $(\theta_\bh, \nu_\bh)$ is the unique solution 
of the linearised system \eqref{eq1_lin}--\eqref{eq4_lin}
in the sense of Theorem~\ref{thm3}.

It remains to show now the strong convergence of $\frac{\varphi_\delta-\varphi}{\delta}$.
To this end, note that by the Lipschitz-contiinuity of $B$ in {\bf A3} 
and \eqref{conv2_lin}
we have 
\[
  \norm{\frac{B(\varphi_\delta)-B(\varphi)}{\delta}}_{L^{p/3}(\Omega; L^\infty(0,T; \cL^2(K,H)))}\leq c\,,
\]
from which, thanks to the classical result \cite[Lem.~2.1]{fland-gat} we get 
\[
  \norm{\int_0^\cdot\frac{B(\varphi_\delta(s))-
  B(\varphi(s))}{\delta}\,\d W(s)}_{L^{p/3}_\cP(\Omega; W^{r,p/3}(0,T; H))} \leq c_r
  \qquad\forall\,r\in(0,1/2)\,.
\]
By comparison in the equation \eqref{var_lin_app}
and the estimates proved above we infer then that 
\[
  \norm{\frac{\varphi_\delta-\varphi}\delta}_{L^{p/3}(\Omega; W^{r,p/3}(0,T; V_2^*))} \leq c_r
  \qquad\forall\,r\in(0,1/2)\,.
\]
Now, recalling that by \cite[Cor.~5]{simon} we have
\[
  L^2(0,T; V_2)\cap W^{r,p/3}(0,T; V_2^*) \embed 
  L^2(0,T; V_1)
  \qquad\text{compactly}\,,
\]
so that the laws of $(\frac{\varphi_\delta-\varphi}\delta)_\delta$
are tight on $L^2(0,T; V_1)$.
By using again Lemma~\ref{strong} together with the uniqueness of 
the limit problem at $\delta=0$, 
proceeding as in Subsection~\ref{ssec:pass_lim}
we also get the strong convergence
\[
  \frac{\varphi_\delta-\varphi}\delta \to \theta_\bh \qquad\text{in } 
  L^2(0,T; V_1)\,,
  \quad\P\text{-a.s.}
\]
which in turn yields, together with 
\eqref{conv2_lin}, the strong convergence of Theorem~\ref{thm3}.
This proves that
$S_1$ is G\^ateaux-differentiable
and its derivative is a solution to the linearised system.
in the sense of Theorem~\ref{thm3}.

\subsection{Fr\'echet-differentiability}
We are only left to show the Fr\'echet-differentiability of $S_1$. To this end, 
since $\widetilde\mU_{ad}$ is open in $\mU$, there is a 
$\mU$-ball $B^\mU_r(\bu)$ of 
radius $r=r_\bu>0$ centred at $\bu$ such that
$B^\mU_r(\bu)\subset\widetilde\mU_{ad}$.
For all $\bh\in B^\mU_r({\bf 0})$, we set $(\varphi_\bh,\mu_\bh):=S(\bu+\bh)$,
$y_\bh:=\varphi_\bh-\varphi-\theta_\bh$, and $z_\bh:=\mu_\bh-\mu-\nu_\bh$, so that 
\begin{align*}
  &\d y_\bh - \Delta z_\bh\,\d t  
  +\bu\cdot\nabla y_\bh\,\d t 
  +\bh\cdot\nabla(\varphi_\bh-\varphi)\,\d t
  = (B(\varphi_\bh)-B(\varphi)-DB(\varphi)\theta_\bh)\,\d W\,,\\
  &z_\bh=-\Delta y_\bh + F'(\varphi_\bh) - F'(\varphi) - F''(\varphi)\theta_\bh\,.
\end{align*}
Noting that $(y_\bh)_\OO=0$, It\^o's formula yields 
\begin{align*}
  &\frac12\norm{\nabla\mN y_\bh(t)}_H^2 + 
  \int_{Q_t}|\nabla y_\bh|^2
  +\int_{Q_t}(F'(\varphi_\bh) - F'(\varphi) - F''(\varphi)\theta_\bh)y_\bh
  -\int_{Q_t}(\varphi_\bh-\varphi)\bh\cdot\nabla\mN y_\bh\\
  &=\int_{Q_t}y_\bh\bu\cdot\nabla\mN y_\bh 
  +\int_0^t\left(\mN y_\bh(s),
  (B(\varphi_\bh(s))-B(\varphi(s))-DB(\varphi(s))\theta_\bh(s))\,\d W(s)\right)_H\\
  &\qquad+\frac12\int_0^t\norm{\nabla\mathcal N(B(\varphi_\bh(s))-
  B(\varphi(s))-DB(\varphi(s))\theta_\bh(s))}_{\cL^2(K,H)}^2\,\d s
  \qquad\forall\,t\in[0,T]\,,\quad\P\text{-a.s.}
\end{align*}
Now, the Young and H\"older inequalities give,
together with the embedding $V_1\embed L^6(\OO)$,
\[
  \int_{Q_t}y_\bh\bu\cdot\nabla\mN y_\bh \leq\eps\int_{Q_t}|\nabla y_\bh|^2
  +c_\eps\int_0^t\left(1+\norm{\bu(s)}_U^2\right)\norm{\nabla\mN y_\bh(s)}_H^2\,\d s
  \qquad\forall\,\eps>0
\]
and similarly
\[
  \int_{Q_t}(\varphi_\bh-\varphi)\bh\cdot\nabla\mN y_\bh\leq
  \int_{Q_t}|\nabla\mN y_\bh|^2 + c\norm{\varphi_\bh-\varphi}_{L^4(0,T; V_1)}^2
  \norm{\bh}_{L^4(0,T; U)}^2\,.
\]
Moreover, note that by the mean value theorem and assumption {\bf A1} we have 
\begin{align*}
  &\int_{Q_t}(F'(\varphi_\bh) - F'(\varphi) - F''(\varphi)\theta_\bh)y_\bh\\
  &=\int_{Q_t}\int_0^1F''(\varphi+\sigma(\varphi_\bh-\varphi))
  |y_\bh|^2\,\d\sigma + 
  \int_{Q_t}\int_0^1\left(F''(\varphi+\sigma(\varphi_\bh-\varphi))-F''(\varphi)\right)
  \theta_\bh y_\bh\,\d\sigma\\
  &\geq-C_\Psi\int_{Q_t}|y_\bh|^2
  +\int_{Q_t}\int_0^1\int_0^1
  F'''(\varphi+\sigma\tau(\varphi_\bh-\varphi))\sigma(\varphi_\bh-\varphi)
  \theta_\bh y_\bh\,\d\tau\,\d\sigma\,,
\end{align*}
where, by the H\"older inequality, the compactness inequality 
\eqref{comp_ineq},
the embedding $V_1\embed L^6(\OO)$, and assumption {\bf C1},
\begin{align*}
  &\int_{Q_t}\int_0^1\int_0^1
  F'''(\varphi+\sigma\tau(\varphi_\bh-\varphi))\sigma(\varphi_\bh-\varphi)
  \theta_\bh y_\bh\,\d\tau\,\d\sigma\\
  &\leq c\int_0^t\left(1+\norm{\varphi(s)}_{L^6(\OO)}
  +\norm{\varphi_\bh(s)}_{L^6(\OO)}\right)
  \norm{(\varphi_\bh-\varphi)(s)}_{L^6(\OO)}
  \norm{\theta_\bh(s)}_{L^6(\OO)}\norm{y_\bh(s)}_H \,\d s
  \\
  &\leq\eps \int_{Q_t}|\nabla y_\bh|^2 + c_\eps\int_{Q_t}|\nabla\mN y_\bh|^2\\
  &\qquad+c\left(1+\norm{\varphi}_{L^\infty(0,T;V_1)}^2
  +\norm{\varphi_\bh}^2_{L^\infty(0,T; V_1)}\right)
  \norm{\varphi-\varphi_\bh}^2_{L^4(0,T; V_1)}
  \norm{\theta_\bh}^2_{L^4(0,T; V_1)}\,.
\end{align*}
Lastly, we have 
\begin{align*}
  B(\varphi_\bh)-B(\varphi)-DB(\varphi)\theta_\bh=
  \int_0^1\left[DB(\varphi+\sigma(\varphi_\bh-\varphi))y_\bh + 
  \left(DB(\varphi+\sigma(\varphi_\bh-\varphi))-DB(\varphi)\right)\theta_\bh\right]\,\d\sigma
\end{align*}
so that by {\bf A3}, {\bf C2--C3}, and the compactness inequality \eqref{comp_ineq},
\begin{align*}
  &\frac12\int_0^t\norm{B(\varphi_\bh(s))-B(\varphi(s))-DB(\varphi(s))
  \theta_\bh(s)}_{\cL^2(K,H)}^2\,\d s\\
  &\leq C_B^2\int_{Q_t}|y_\bh|^2+
   c\int_0^t\norm{(\varphi_\bh-\varphi)(s)}_{V_1}^2\norm{\theta_\bh(s)}^2_{V_1}\,\d s\\
   &\leq\eps\int_{Q_t}|\nabla y_\bh|^2 + c_\eps\int_{Q_t}|\nabla\mN y_\bh|^2
   +c\norm{\varphi-\varphi_\bh}^2_{L^4(0,T; V_1)}\norm{\theta_\bh}^2_{L^4(0,T; V_1)}\,.
\end{align*}
Consequently, taking all this information into account, 
we can choose $\eps$ small enough and rearrange the terms to get 
\begin{align*}
  &\frac12\norm{\nabla\mN y_\bh(t)}^2_H + \int_{Q_t}|\nabla y_\bh|^2\\
  &\leq
  \int_0^t\left(1+\norm{\bu(s)}_U^2\right)\norm{\nabla\mN y_\bh(s)}_H^2\,\d s
  +c\norm{\varphi_\bh-\varphi}_{L^4(0,T; V_1)}^2\norm{\bh}_{L^4(0,T; U)}^2\\
  &\qquad+
  c\left(1+\norm{\varphi}_{L^\infty(0,T;V_1)}^2
  +\norm{\varphi_\bh}^2_{L^\infty(0,T; V_1)}\right)
  \norm{\varphi-\varphi_\bh}^2_{L^4(0,T; V_1)}
  \norm{\theta_\bh}^2_{L^4(0,T; V_1)}\\
  &\qquad+\int_0^t\left(\mN y_\bh(s),
  (B(\varphi_\bh(s))-B(\varphi(s))-DB(\varphi(s))\theta_\bh(s))\,\d W(s)\right)_H\,.
\end{align*}
Thanks to the embedding 
$L^\infty(0,T; H)\cap L^2(0,T; V_2)\embed L^4(0,T; V_1)$,
by \eqref{cont2} and \eqref{conv1_lin}--\eqref{conv2_lin} we have 
\[
  \norm{\varphi_\bh-\varphi}_{L^{p/3}_\cP(\Omega; L^4(0,T; V_1))}
  +\norm{\theta_\bh}_{L^{p/3}_\cP(\Omega; L^4(0,T; V_1))}\leq c\norm{\bh}_{\mU}\,,
\]
while \eqref{bound} yields
\[
  \norm{\varphi_\bh}_{L^{p}_\cP(\Omega; L^\infty(0,T; V_1))}
  +\norm{\varphi}_{L^{p}_\cP(\Omega; L^\infty(0,T; V_1))}\leq c\,,
\]
where the constant $c$ is independent of $\bh$. 
Taking power $\frac{p}{14}$ at both sides, supremum in time and expectations, 
on the right-hand side we use the H\"older inequality with exponents 
$\frac17+\frac37+\frac37=1$ to get 
\[
\norm{\left(1+\norm{\varphi}_{L^\infty(0,T;V_1)}^{p/7}
  +\norm{\varphi_\bh}^{p/7}_{L^\infty(0,T; V_1)}\right)
  \norm{\varphi-\varphi_\bh}^{p/7}_{L^4(0,T; V_1)}
  \norm{\theta_\bh}^{p/7}_{L^4(0,T; V_1)}}_{L^{1}(\Omega)}\leq
  c\norm{\bh}_{\mU}^{2p/7}
\]
and similarly
\[
  \norm{\norm{\varphi_\bh-\varphi}_{L^4(0,T; V_1)}^{p/7}
  \norm{\bh}_{L^4(0,T; U)}^{p/7}}_{L^1(\Omega)}\leq
  c\norm{\bh}_{\mU}^{2p/7}\,.
\]
Consequently, arguing again as in Subsection~\ref{ssec:uniq}
using an iterative argument and the Burkholder-Davis-Gundy and Young inequalities
(see also \cite[Lem.~4.1]{mar-scar-ref}) gives then
\[
  \norm{y_\bh}_{L^{p/7}(\Omega; C^0([0,T]; V_1^*)\cap L^2(0,T; V_1))}
  \leq c\norm{\bh}_{\mU}^2=o\left(\norm{\bh}_{\mU}\right)
  \qquad\text{as } \norm{\bh}_\mU\to0\,.
\]
This proves the Fr\'echet-differentiability of $S_1$ and concludes the proof of Theorem~\ref{thm3}.

%%%%%%%%%%%%%%%%%%%%%%%%%%%%%%%%%%

\section{Adjoint system}
\label{sec:ad}
In this section we study the adjoint problem \eqref{eq1_ad}--\eqref{eq4_ad},
proving that it is well-posed
in the sense of Theorem~\ref{thm4}.

As we have anticipated in the introduction, the presence of the extra--random
component in the convection term calls for non--trivial 
mathematical tools when deriving estimates on the solutions.
Let us recall here a general backward version of the 
stochastic Grownall lemma that will be used in this section: for details
we refer to 
\cite[Thm.~1]{hun-gronwall} and \cite{wang-fan}.
\begin{lem}
 \label{lem:gronwall}
 Let $\xi\in L^2(\Omega,\cF_T)$ be non-negative, 
 $\alpha\in L^\infty_\cP(\Omega; L^1(0,T))$ 
 with $\alpha\geq\alpha_0>0$
 almost everywhere in $\Omega\times(0,T)$,
 and $X\in L^2_\cP(\Omega; C^0([0,T]))$
 be a non-negative process such that 
 \[
 X(t) \leq \E\left[\xi + \int_t^T\alpha(s)X(s)\,\d s\,\Bigg|\cF_t\right]
 \qquad\forall\,t\in[0,T]\,,\quad\P\text{-a.s.}
 \]
 Then, for every $t\in[0,T]$ it holds that 
 \[
 X(t) \leq \E\left[\xi\exp\norm{\alpha}_{L^1(t,T)}\,\Big|\cF_t\right] \qquad\P\text{-a.s.}
 \]
\end{lem}

\subsection{Approximation}
For every $\lambda>0$, using the approximations on $\Psi$ and $\bu$
as in Section~\ref{ssec:approx}, we consider the approximated problem 
\begin{align}
  \nonumber
  -\d P_\lambda -\Delta \tilde P_\lambda \,\d t 
  + \Psi_\lambda''(\varphi)\tilde P_\lambda\,\d t - \bu_\lambda\cdot\nabla P_\lambda\,\d t
  \qquad\qquad& \\
  \label{eq1_ad_lam}
  =\alpha_1(\varphi-\varphi_Q)\,\d t+
  DB(\varphi)^*Z_\lambda\,\d t  - Z_\lambda\,\d W
  \qquad&\text{in } (0,T)\times\OO\,,\\
  \label{eq2_ad_lam}
  \tilde P_\lambda=-\Delta P_\lambda
  \qquad&\text{in } (0,T)\times\OO\,,\\
  \label{eq3_ad_lam}
  {\bf n}\cdot\nabla P_\lambda = {\bf n}\cdot\nabla \tilde P_\lambda = 0
  \qquad&\text{in } (0,T)\times\partial\OO\,,\\
  \label{eq4_ad_lam}
  P_\lambda(T)=\alpha_2(\varphi(T)-\varphi_T)
  \qquad&\text{in } \OO\,.
\end{align}
This can be written in abstract form as
\[
  -\d P_\lambda + \mathcal F_\lambda(P_\lambda)\,\d t=
  \alpha_1(\varphi-\varphi_Q)\,\d t + 
   DB(\varphi)^*Z_\lambda\,\d t - Z_\lambda\,\d W\,, \qquad
  P_\lambda(T)=\alpha_2(\varphi(T)-\varphi_T)\,,
\]
where 
$\mathcal F_\lambda:\Omega\times[0,T]\times V_2\to V_2^*$
is given by 
\[
  \ip{\mathcal F_\lambda(\omega,t,y)}{\zeta}:=
  \int_\OO\left(\Delta y\Delta\zeta
  -\Psi_\lambda''(\varphi(\omega,t))\Delta y\zeta + 
  y\bu_\lambda(\omega,t)\cdot\nabla\zeta\right)\,,\qquad y,\zeta\in V_2\,.
\]
By construction it holds that 
$\Psi''_\lambda(\varphi)\in L^\infty(\Omega\times Q)$
and $\bu_\lambda\in L^\infty_\cP(\Omega\times(0,T); U)$, so that 
using similar arguments to the ones in Subsection~\ref{ssec:approx}, we have that
the operator $\mathcal F_\lambda$ is progressively measurable, 
hemicontinuous, weakly  monotone, 
weakly coercive, and linearly bounded.
Moreover, the Lipschitz-continuity of $B$ in {\bf A3} implies that 
$DB(\varphi)^*$ is uniformly bounded as well.
The classical variational theory for backward SPDEs \cite[Sec.~3]{du-meng2}
ensures then that such approximated problem admits a unique 
variational solution $(P_\lambda, Q_\lambda)$, with 
\begin{align*}
  P_\lambda 
  \in L^2_{\cP}(\Omega; C^0([0,T]; H)\cap  L^2(0,T; V_2))\,, \qquad
  Z_\lambda \in L^2_\cP(\Omega; L^2(0,T; \cL^2(U,H)))\,.
\end{align*}
Actually, let us note that 
thanks to the assumption on the target $\varphi_T$
and the regularity of $\varphi$,
the final value satisfies $\alpha_2(\varphi(T)-\varphi_T)\in L^2(\Omega, \cF_T; V_1)$.
Consequently, by a standard finite dimensional approximation 
of the approximated problem with $\lambda>0$ fixed, it follows that 
the approximated solution actually inherits more regularity, namely
\begin{align*}
  P_\lambda 
  \in L^2_{\cP}(\Omega; C^0([0,T]; V_1)\cap  L^2(0,T; V_3))\,, \qquad
  Z_\lambda \in L^2_\cP(\Omega; L^2(0,T; \cL^2(U,V_1)))\,.
\end{align*}
We can then set
\[
  \tilde P_\lambda :=\mathcal LP_\lambda 
  \in L^2_{\cP}(\Omega; C^0([0,T]; V_1^*)\cap  L^2(0,T; V_1))\,,
\]
so that $(P_\lambda, \tilde P_\lambda, Z_\lambda)$ satisfy,
for every $t\in[0,T]$, $\P$-almost surely, for every $\zeta\in V_1$,
\begin{align*}
    &\left(P_\lambda(t), \zeta\right)_H
    +\int_{Q_t^T}\nabla\tilde P_\lambda\cdot\nabla\zeta
    +\int_{Q_t^T}\Psi_\lambda''(\varphi)\tilde P_\lambda\zeta
    +\int_{Q_t^T}P_\lambda\bu_\lambda\cdot\nabla\zeta\\
    &\quad
    =\left(\alpha_2(\varphi(T)-\varphi_T), \zeta\right)_H
    +\int_{Q_t^T}\alpha_1(\varphi-\varphi_Q)\zeta
    +\int_{Q_t^T}DB(\varphi)^*Z_\lambda\zeta
    -\left(\int_t^TZ_\lambda(s)\,\d W(s), \zeta\right)_H\,.
\end{align*}

\subsection{An estimate by duality method}
The first estimate that we prove is based on a duality method 
between the approximated adjoint system \eqref{eq1_ad_lam}--\eqref{eq4_ad_lam}
and a suitably introduced approximated linearised system.
This step is fundamental as it allows to obtain some 
preliminary estimates on the adjoint variables
without working explicitly on the adjoint system, which may be 
not trivial. Such duality method is extremely powerful,
and it will be crucial in showing well-posedness of the adjoint system.

The idea is the following: we consider 
the $\lambda$-approximated version 
of the linearised system \eqref{eq1_lin}--\eqref{eq4_lin},
in a more general version where the forcing term is given by 
an arbitrary term 
\[
g\in L^{\frac{2p}{p+4}}_\cP(\Omega; L^2(0,T; H))\,.
\]
Namely, for $\bh\in\mathcal U$ we consider 
\begin{align}
    \label{eq1_lin_lam}
    \d\theta_{\bh,\lambda}^g - \Delta \nu_{\bh,\lambda}^g\,\d t 
    +\bh\cdot\nabla\varphi\,\d t
    + \bu_\lambda\cdot\nabla\theta^g_{\bh,\lambda}\,\d t
    = DB(\varphi)\theta^g_{\bh,\lambda}\,\d W \qquad&\text{in } (0,T)\times \OO\,,\\
    \label{eq2_lin_lam}
    \nu_{\bh,\lambda}^g=-\Delta \theta_{\bh,\lambda}^g +
     \Psi_\lambda''(\varphi)\theta_{\bh,\lambda}^g -g
     \qquad&\text{in } (0,T)\times \OO\,,\\
    \label{eq3_lin_lam}
    {\bf n}\cdot\nabla\theta_{\bh,\lambda}^g = {\bf n}\cdot\nabla\nu_{\bh,\lambda}^g = 0 
    \qquad&\text{in } (0,T)\times\partial \OO\,,\\
    \label{eq4_lin_lam}
     \theta_{\bh,\lambda}^g(0)=0 \qquad&\text{in } \OO\,.
\end{align}
Since $\Psi''_\lambda(\varphi)\in L^\infty(\Omega\times Q)$, 
the classical variational approach (see again Subsections~\ref{ssec:approx} 
and \ref{ssec:ex_lin}) ensures that
the system \eqref{eq1_lin_lam}--\eqref{eq4_lin_lam}
admits a unique solution
\begin{align*}
  &\theta_{\bh,\lambda}^g \in 
  L^{\frac{2p}{p+4}}_{\cP}\left(\Omega; C^0([0,T]; H)\cap L^2(0,T; V_2)\right)\,,\qquad
  \nu_{\bh,\lambda}^g\in L^{\frac{2p}{p+4}}_{\cP}(\Omega; L^2(0,T; H))\,.
\end{align*}
Moreover, we can show that the system \eqref{eq1_lin_lam}--\eqref{eq4_lin_lam}
is in duality with the approximated adjoint system \eqref{eq1_ad_lam}--\eqref{eq4_ad_lam}.
To this end, by It\^o's formula we have that 
\begin{align*}
  \d(\theta_{\bh,\lambda}^g,P_\lambda)_H&=
  -\tilde P_\lambda \nu_{\bh,\lambda}^g\,\d t
  +\varphi\bh\cdot\nabla P_\lambda\,\d t
  +\theta_{\bh,\lambda}^g\bu_\lambda\cdot\nabla P_\lambda\,\d t
  +(P_\lambda, DB(\varphi)\theta_{\bh,\lambda}^g\,\d W)_H\\
  &+\tilde P_\lambda(-\Delta\theta_{\bh,\lambda}^g+
  \Psi''_\lambda(\varphi)\theta_{\bh,\lambda}^g)\,\d t
  +P_\lambda\bu_\lambda\cdot\nabla\theta_{\bh,\lambda}^g\,\d t
  -\alpha_1(\varphi-\varphi_Q)\theta_{\bh,\lambda}^g\,\d t\\
  &-(DB(\varphi)^*Z_\lambda, \theta_{\bh,\lambda}^g)_H\,\d t 
  + (\theta_{\bh,\lambda}^g, Z_\lambda\,\d W)_H
  +(DB(\varphi)\theta_{\bh,\lambda}^g, Z_\lambda)_{\cL^2(K,H)}\,\d t\,,
\end{align*}
which readily implies by comparison in the two systems that 
\beq
  \label{duality_lam}
  \alpha_1\E\int_Q\theta_{\bh,\lambda}^g(\varphi-\varphi_Q)+
  \alpha_2\E\int_\OO\theta_{\bh,\lambda}^g(T)(\varphi(T)-\varphi_T)=
  \E\int_Q\varphi\bh\cdot\nabla P_\lambda+
  \E\int_Q\tilde P_\lambda g\,.
\eeq

Let us set now for brevity of notation 
$\theta_\lambda^g:=\theta^g_{\bh,\lambda}$ and $\nu_\theta^g:=\nu_{\bh,\lambda}^g$
with the choice $\bh=0$.
Noting that $(\theta_\lambda^g)_\OO=0$, It\^o's formula for 
$\frac12\norm{\nabla\mN\theta_\lambda^g}_H^2$ yields
\begin{align*}
  &\frac12\norm{\nabla\mN\theta_\lambda^g(t)}_H^2 + 
  \int_{Q_t}|\nabla\theta_\lambda^g|^2 = 
  \int_{Q_t}\theta_\lambda^g\bu_\lambda\cdot\nabla\mN\theta_\lambda^g
  -\int_{Q_t}\Psi''_\lambda(\varphi)|\theta_\lambda^g|^2
  +\int_{Q_t}g\theta_\lambda^g\\
  &\qquad+\frac12\int_0^t\norm{\nabla\mN DB(\varphi(s))\theta_\lambda^g(s)}_{\cL^2(K,H)}^2\,\d s
  +\int_0^t\left(\mN\theta_\lambda^g(s),DB(\varphi(s))\theta_\lambda^g(s)\,\d W(s)\right)_H\,.
\end{align*}
Using the fact that $\Psi''_\lambda\geq-C_\Psi$ 
and the boundedness of $DB(\varphi)$ in $\cL(V_1, \cL^2(K,H))$,
thanks to the H\"older--Young inequalities and the compactness 
inequality \eqref{comp_ineq} we get, for all $\eps>0$,
\begin{align*}
  \norm{\nabla\mN\theta_\lambda^g(t)}_H^2 + 
  \int_{Q_t}|\nabla\theta_\lambda^g|^2 &\leq 
  \int_Q|g|^2 + \eps\int_{Q_t}|\nabla\theta_\lambda^g|^2
  + c_\eps
  \int_0^t\left(1+\norm{\bu(s)}_U^2\right)\norm{\nabla\mN\theta_\lambda^g(s)}_H^2\,\d s\\
  &+\int_0^t\left(\mN\theta_\lambda^g(s),DB(\varphi(s))\theta_\lambda^g(s)\,\d W(s)\right)_H\,.
\end{align*}
We take now power $\frac{p}{p+4}$ at both sides, supremum in time, and expectations.
Thanks to the Burkholder--Davis--Gundy inequality
(see \cite[Lem.~4.1]{mar-scar-ref}), assumption {\bf C2},
and \eqref{comp_ineq} we get
\begin{align*}
  \E\sup_{r\in[0,t]}\left|\int_0^r\left(\mN\theta_\lambda^g(s),
  DB(\varphi(s))\theta_\lambda^g(s)\,\d W(s)\right)_H\right|^{\frac{p}{p+4}}
  \leq\frac12\E\norm{\nabla\mN\theta_\lambda^g}_{L^\infty(0,t; H)}^{\frac{2p}{p+4}}
  +c\E\norm{\theta_\lambda^g}_{L^2(0,t; H)}^{\frac{2p}{p+4}}&\\
  \leq\frac12\E\norm{\nabla\mN\theta_\lambda^g}_{L^\infty(0,t; H)}^{\frac{2p}{p+4}}
  +\frac12\E\norm{\nabla\theta_\lambda^g}_{L^2(0,t; H)}^{\frac{2p}{p+4}}
  +ct^{\frac{p}{p+4}}\E\norm{\nabla\mN\theta_\lambda^g}_{L^\infty(0,t; H)}^{\frac{2p}{p+4}}&\,.
\end{align*}
Moreover, 
since $\bu\in\widetilde\mU_{ad}$, by the H\"older inequality we have
\begin{align*}
  &\E\sup_{r\in[0,t]}\left|\int_0^r\left(1+\norm{\bu(s)}_U^2\right)
  \norm{\nabla\mN\theta_\lambda^g(s)}_H^2\,\d s\right|^{\frac{p}{p+4}}\\
  &\leq c\E\left|
  t^{1-\frac2p}\left(1+\norm{\bu}_{L^p(0,T; U)}^2\right)
  \norm{\nabla\mN\theta_{\lambda}^g}_{L^\infty(0,t; H)}^2\right|^{\frac{p}{p+4}}
  \leq
  ct^{\frac{p-2}{p+4}}\E\norm{\nabla\mN
  \theta_{\lambda}^g}_{L^\infty(0,t; H)}^{\frac{2p}{p+4}}\,.
\end{align*}
Since $\frac{p}{p+4}>0$ and $\frac{p-2}{p+4}>0$, 
we can close the estimate rearranging all the terms on $[0,T_0]$
for $T_0$ sufficiently small (independent of both $\lambda$ and $g$).
Using once more a classical iterative procedure on 
every subinterval until $T$, we infer that
there exists a constant $c>0$, independent of both $\lambda$ and $g$, 
such that 
\beq\label{est1_ad}
  \norm{\theta_\lambda^g}_{L^{\frac{2p}{p+4}}_\cP(\Omega; 
  C^0([0,T]; V_1^*)\cap L^2(0,T; V_1))} 
  \leq c\norm{g}_{L^{\frac{2p}{p+4}}_\cP(\Omega; L^2(0,T; H))}\,.
\eeq
Now,
by assumption {\bf C4} and the regularity of $\varphi$
(since $\frac{2p}{p-4}\leq p$ for $p\geq6$) it holds 
\[
\alpha_1(\varphi-\varphi_Q)\in L^{\frac{2p}{p-4}}_\cP(\Omega; L^2(0,T; H))\,,\qquad
\alpha_2(\varphi(T)-\varphi_T)\in L^{\frac{2p}{p-4}}(\Omega,\cF_T; V_1)\,,
\]
so that the duality relation \eqref{duality_lam}
(with $\bh=0$) and the estimate \eqref{est1_ad}
yield
\begin{align*}
  \E\int_Q\tilde P_\lambda g&\leq
  \norm{\theta_\lambda^g}_{L^{\frac{2p}{p+4}}_\cP(\Omega; L^2(0,T; H))}
  \norm{\alpha_1(\varphi-\varphi_Q)}_{L^{\frac{2p}{p-4}}_\cP(\Omega; L^2(0,T; H))} \\
  &+
  \norm{\theta_\lambda^g}_{L^{\frac{2p}{p+4}}_\cP(\Omega; 
  C^0([0,T]; V_1^*))}
  \norm{\alpha_2(\varphi(T)-\varphi_T)}_{L^{\frac{2p}{p-4}}(\Omega,\cF_T;V_1)}
  \leq c\norm{g}_{L^{\frac{2p}{p+4}}_\cP(\Omega; L^2(0,T; H))}\,.
\end{align*}
By the arbitrariness of $g$ we obtain 
\beq
  \label{est2_ad}
  \|\tilde P_\lambda\|_{L^{\frac{2p}{p-4}}_\cP(\Omega; L^2(0,T; H))} \leq c\,.
\eeq

\subsection{Further estimates}
\label{ssec:further}
We show here that the initial estimate \eqref{est2_ad} allows to 
obtain uniform estimates on the adjoint variables.
To this end, It\^o's formula for $\frac12\norm{P_\lambda}_H^2+
\frac12\norm{\nabla P_\lambda}_H^2$ yields,
recalling that $\tilde P_\lambda=\mathcal L P_\lambda$,
\begin{align}
  \nonumber
  &\frac12\norm{P_\lambda(t)}_{V_1}^2
  +\int_t^T\|\tilde P_\lambda(s)\|_{V_1}^2\,\d s 
  +\frac12\int_t^T\norm{Z_\lambda(s)}^2_{\cL^2(K,V_1)}\,\d s\\
  \nonumber
  &=\frac{\alpha_2^2}2\norm{\varphi(T)-\varphi_T}_{V_1}^2
  -\int_{Q_t^T}\Psi''_\lambda(\varphi)|\tilde P_\lambda|^2
  -\int_{Q_t^T}\Psi''_\lambda(\varphi)\tilde P_\lambda P_\lambda
  +\int_{Q_t^T}(P_\lambda+ \tilde P_\lambda)\bu_\lambda\cdot\nabla P_\lambda\\
  \nonumber
  &\qquad+\alpha_1\int_{Q_t^T}(\varphi-\varphi_Q)(P_\lambda+ \tilde P_\lambda)
  +\int_{Q_t^T}(DB(\varphi)^*Z_\lambda)(P_\lambda+ \tilde P_\lambda)\\
  &\qquad
  -\int_t^T\left(P_\lambda(s) + \tilde P_\lambda(s), Z_\lambda(s)\,\d W(s)\right)_H
  \qquad\forall\,t\in[0,T]\,,\quad\P\text{-a.s.}
  \label{aux1_ad}
\end{align}
On the right-hand side, we have already noticed that 
$\alpha_2(\varphi(T)-\varphi_T)\in L^2(\Omega,\cF_T; V_1)$.
Moreover, by {\bf A1}, the compactness inequality \eqref{comp_ineq}
and the fact that $\tilde P_\lambda=\mathcal L P_\lambda$, 
for the second and third terms we have
\[
  -\int_{Q_t^T}\Psi''_\lambda(\varphi)|\tilde P_\lambda|^2 \leq 
  C_\Psi \int_{Q_t^T}|\tilde P_\lambda|^2\leq
  \eps\int_{Q_t^T}|\nabla\tilde P_\lambda|^2 + c_\eps\int_{Q_t}|\nabla P_\lambda|^2
\]
and, thanks to the H\"older-Young inequalities, the embedding $V_1\embed L^6(\OO)$,
and {\bf C1},
\begin{align*}
  -\int_{Q_t^T}\Psi''_\lambda(\varphi)\tilde P_\lambda P_\lambda
  &\leq\int_t^T\norm{P_\lambda(s)}_{V_1}^2\,\d s + 
  c\int_t^T\norm{\Psi''_\lambda(\varphi(s))}_{L^3(\OO)}^2\|\tilde P_\lambda(s)\|_H^2\,\d s\\
  &\leq\int_t^T\norm{P_\lambda(s)}_{V_1}^2\,\d s +
  c\left(1+\norm{\varphi}_{L^\infty(0,T; V_1)}^4\right)\|\tilde P_\lambda\|^2_{L^2(0,T; H)}\,.
\end{align*}
Also, 
note that since $\tilde P_\lambda=\mathcal L P_\lambda$, 
in particular it holds that $(\tilde P_\lambda)_\OO=0$.
Hence, using the Young and H\"older inequalities,
the embedding $V_1\embed L^6(\OO)$ yields, for all $\eps>0$,
\[
  \int_{Q_t^T}(P_\lambda+\tilde P_\lambda)\bu_\lambda\cdot\nabla P_\lambda
  \leq\eps\int_t^T\|\tilde P_\lambda(s)\|^2_{V_1}\,\d s
  +c_\eps\int_t^T\left(1+\norm{\bu(s)}_U^2\right)\norm{P_\lambda(s)}^2_{V_1}\,\d s\,,
\]
and similarly 
\[
  \alpha_1\int_{Q_t^T}(\varphi-\varphi_Q)(P_\lambda+\tilde P_\lambda) \leq
  \alpha_1^2\int_Q|\varphi-\varphi_Q|^2 + 
  \frac12\int_{Q_t^T}|P_\lambda|^2+
  \frac12\int_{Q_t^T}|\tilde P_\lambda|^2\,.
\]
Lastly, thanks to {\bf A3} and {\bf C2}, and again the compactness inequality \eqref{comp_ineq},
we have that
\begin{align*}
  \int_{Q_t^T}(DB(\varphi)^*Z_\lambda)(P_\lambda+\tilde P_\lambda) &= 
  \int_t^T\left(Z_\lambda(s), DB(\varphi(s))
  (P_\lambda+\tilde P_\lambda)(s)\right)_{\cL^2(K,H)}\,\d s\\
  &\leq \frac14\int_t^T\norm{Z_\lambda(s)}^2_{\cL^2(K,H)}\,\d s
  +2C_B^2\int_{Q_t^T}|P_\lambda|^2+2C_B^2\int_{Q_t^T}|\tilde P_\lambda|^2\\
  &\leq \frac14\int_t^T\norm{Z_\lambda(s)}^2_{\cL^2(K,H)}\,\d s
  +\eps\int_{Q_t^T}|\nabla \tilde P_\lambda|^2+c_\eps\int_t^T\norm{P_\lambda(s)}_{V_1}^2\,\d s\,.
\end{align*}
Choosing $\eps$ small enough, rearranging the terms in \eqref{aux1_ad},
and conditioning \eqref{aux1_ad} with respect to $\cF_t$
we are left with 
\begin{align*}
  &\norm{P_\lambda(t)}_{V_1}^2
  +\E\left[\int_t^T\|\tilde P_\lambda(s)\|_{V_1}^2\,\d s 
  +\int_t^T\norm{Z_\lambda(s)}^2_{\cL^2(K,V_1)}\,\d s\,\Bigg| \cF_t\right]\\
  &\leq c +  c
  \E\left[\left(1+\norm{\varphi}_{L^\infty(0,T; V_1)}^4\right)\|\tilde P_\lambda\|^2_{L^2(0,T; H)}
  +\int_t^T\left(1+\norm{\bu(s)}_U^2\right)\norm{P_\lambda(s)}^2_{V_1}\,\d s\,
  \Bigg| \cF_t\right]\,,
\end{align*}
so that the backward version of the stochastic Gronwall Lemma~\ref{lem:gronwall}
yields
\begin{align*}
  &\norm{P_\lambda(t)}_{V_1}^2
  +\E\left[\int_t^T\|\tilde P_\lambda(s)\|_{V_1}^2\,\d s 
  +\int_t^T\norm{Z_\lambda(s)}^2_{\cL^2(K,V_1)}\,\d s\,\Bigg| \cF_t\right]\\
  &\leq\E\left[\left(c+c\left(1 + \norm{\varphi}_{L^\infty(0,T; V_1)}^4\right)
 \|\tilde P_\lambda\|^2_{L^2(0,T; H)}\right)
 \exp\left(t+\norm{\bu}_{L^2(0,T; U)}^2\right)
  \,\Bigg| \cF_t\right]\,.
\end{align*}
Consequently, taking expectations
we infer that
\begin{align*}
  &\E\norm{P_\lambda(t)}_{V_1}^2 + \E\|\tilde P_\lambda\|^2_{L^2(t,T; V_1)}
  +\E\norm{Z_\lambda}^2_{L^2(t,T; \cL^2(K,V_1))}\\
  &\leq c\left(1+\exp\norm{\bu}_\mU^2\right)
  \E\left[1+\left(1+\norm{\varphi}_{L^\infty(0,T; V_1)}^4\right)
  \|\tilde P_\lambda\|^2_{L^2(0,T; H)}\right]\,,
\end{align*}
where, by the H\"older inequality and the duality-estimate \eqref{est2_ad}, we have 
\begin{align*}
  \E\left[\left(1+\norm{\varphi}_{L^\infty(0,T; V_1)}^4\right)\|\tilde P_\lambda\|^2_{L^2(0,T; H)}\right]
  &\leq\norm{1+\norm{\varphi}_{L^\infty(0,T; V_1)}^4}_{L^{\frac{p}4}(\Omega)}
  \norm{\|\tilde P_\lambda\|^2_{L^2(0,T; H)}}_{L^{\frac{p}{p-4}}(\Omega)}\\
  &\leq c\left(1+\norm{\varphi}^4_{L^p(\Omega; L^\infty(0,T; V_1))}\right)
  \|\tilde P_\lambda\|^2_{L^\frac{2p}{p-4}(\Omega; L^2(0,T; H))}\leq c\,,
\end{align*}
which yields in turn
\beq\label{est3_ad}
  \norm{P_\lambda}_{C^0([0,T]; L^2(\Omega; V_1))} + 
  \|\tilde P_\lambda\|_{L^2_\cP(\Omega; L^2(0,T; V_1))}
  +\norm{Z_\lambda}_{L^2_\cP(\Omega; L^2(0,T; \cL^2(K,V_1)))} \leq c\,.
\eeq
With this additional information, we can perform a classical refinement on the estimates
going back to the inequality \eqref{aux1_ad},
repeating the same steps but this time taking first supremum in time and then expectations:
the estimate \eqref{est3_ad} allows to apply the Burkholder-Davis-Gundy inequality on the 
stochastic integral, so that we obtain, thanks also to elliptic regularity,
\beq
  \label{est4_ad}
  \norm{P_\lambda}_{L^2_\cP(\Omega; C^0([0,T]; V_1)\cap L^2(0,T; V_3))}+
  \|\tilde P_\lambda\|_{L^2_\cP(\Omega; C^0([0,T]; V_1^*)\cap L^2(0,T; V_1))}\leq c\,.
\eeq

\subsection{\bf Passage to the limit}
From \eqref{est3_ad}--\eqref{est4_ad} we infer that
there exists $(P,\tilde P,Z)$ with 
\begin{align*}
  &P\in L^2_w(\Omega; L^\infty(0,T; V_1))\cap L^2_\cP(\Omega; L^2(0,T; V_3))\,,\\
  &\tilde P=\mathcal LP \in 
  L^2_w(\Omega; L^\infty(0,T; V_1^*))\cap L^2_\cP(\Omega; L^2(0,T; V_1))\,,\\
  &Z \in L^2_\cP(\Omega; L^2(0,T; \cL^2(K,V_1)))\,,
\end{align*}
such that as $\lambda\searrow0$, possibly on a subsequence,
\begin{align}
  \label{conv1_ad}
  P_\lambda \wstarto P \qquad&\text{in } 
  L^2_w(\Omega; L^\infty(0,T; V_1))\cap L^2_\cP(\Omega; L^2(0,T; V_3))\,,\\
  \label{conv2_ad}
  \tilde P_\lambda \wstarto \tilde P\qquad&\text{in } 
  L^2_w(\Omega; L^\infty(0,T; V_1^*))\cap L^2_\cP(\Omega; L^2(0,T; V_1))\,,\\
  \label{conv3_ad}
  Z_\lambda \wto Z \qquad&\text{in } L^2_\cP(\Omega; L^2(0,T; \cL^2(K,V_1)))\,. 
\end{align}
Now, thanks to {\bf C1} and the regularity of $\varphi$ we have 
$\Psi''(\varphi) \in L^3(\Omega; L^\infty(0,T; L^3(\OO)))$, so 
in particular
\[
  \Psi''_\lambda(\varphi) \to \Psi''(\varphi) \qquad\text{in } L^3(\Omega\times Q)\,,
\]
and also, thanks to \eqref{conv2_ad},
\[
  \Psi''_\lambda(\varphi)\tilde P_\lambda 
  \wto \Psi''(\varphi)\tilde P \qquad\text{in } L^{6/5}_\cP(\Omega; L^{6/5}(0,T; L^{6/5}(\OO)))\,.
\]
Similarly, since $\bu_\lambda\to\bu$ in $L^q_\cP(\Omega;L^p(0,T; U))$
for every $q\geq1$, from \eqref{conv1_ad} we have 
\[
  \bu_\lambda\cdot\nabla P_\lambda \wto \bu\cdot\nabla P
  \qquad\text{in } L^{\ell}_\cP(\Omega; L^{\frac{2p}{p+2}}(0,T;H))
  \qquad\forall\,\ell\in[1,2)\,. 
\]
Lastly, convergence \eqref{conv3_ad} readily implies that 
\[
  DB(\varphi)^*Z_\lambda \wto DB(\varphi)^*Z \qquad\text{in }
  L^2_\cP(\Omega; L^2(0,T; H))\,,
\]
while by the linearity and continuity of
the stochastic integral we have 
\[
  \int_\cdot^TZ_\lambda(s)\,\d W(s) \wto 
  \int_\cdot^TZ(s)\,\d W(s) \qquad\text{in } L^2_\cP(\Omega; C^0([0,T]; V_1))\,.
\]
Consequently, we can let $\lambda\searrow0$ in the 
variational formulation of the approximated system 
\eqref{eq1_ad_lam}--\eqref{eq4_ad_lam} and deduce that $(P,\tilde P, Z)$
solve the limit adjoint problem \eqref{eq1_ad}--\eqref{eq4_ad}.
The pathwise continuity of $P$, hence by comparison also of $\tilde P$,
follows by classical methods using It\^o's formula on the limit equation.

\subsection{Uniqueness}
By linearity of the adjoint system, 
it is enough to show that if $(P,\tilde P, Z)$ is a solution of 
\eqref{eq1_ad}--\eqref{eq4_ad} with $\alpha_1=\alpha_2=0$,
then $\nabla P=0$, $\tilde P=0$, and $\nabla Z=0$.
To this end, It\^o's formula for $\frac12\norm{\nabla P}_H^2$ yields
\begin{align*}
  \nonumber
  &\frac12\norm{\nabla P(t)}_{H}^2
  +\int_{Q_t^T}|\nabla\tilde P|^2
  +\frac12\int_t^T\norm{\nabla Z(s)}^2_{\cL^2(K,H)}\,\d s\\
  &=
  -\int_{Q_t^T}\Psi''(\varphi)|\tilde P|^2
  +\int_{Q_t^T}\tilde P \bu\cdot\nabla P
  +\int_{Q_t^T}(DB(\varphi)^*Z)\tilde P
  -\int_t^T\left(\tilde P(s), Z(s)\,\d W(s)\right)_H
\end{align*}
Now, as the computations are similar to the ones of Subsection~\ref{ssec:further},
we avoid details for brevity.
The terms on the right-hand side can be treated 
using {\bf A1}, the H\"older--Young inequalities, the embedding $V_1\embed L^6(\OO)$,
and the compactness inequality \eqref{comp_ineq} as
\[
  -\int_{Q_t^T}\Psi''(\varphi)|\tilde P|^2
  +\int_{Q_t^T}\tilde P \bu\cdot\nabla P\leq \eps\int_{Q_t^T}|\nabla\tilde P|^2
  +c_\eps\int_0^t\left(1+\norm{\bu(s)}_U^2\right)\norm{\nabla P(s)}_H^2\,\d s\,,
\]
and similarly, since $DB(\varphi)\tilde P$ is $\cL^2(K,H_0)$-valued by {\bf A3},
by the Poincar\'e--Wirtinger inequality and {\bf C2} we have
\begin{align*}
  \int_{Q_t^T}(DB(\varphi)^*Z)\tilde P&=
  \int_0^t(Z(s), DB(\varphi(s))\tilde P(s))_{\cL^2(K,H)}\,\d s\\
 &\leq \frac14\int_t^T\norm{\nabla Z(s)}^2_{\cL^2(K,H)}\,\d s
 +\eps\int_{Q_t^T}|\nabla\tilde P|^2 + c_\eps\int_t^T\norm{\nabla P(s)}_H^2\,\d s\,.
\end{align*}
Rearranging the terms and taking conditional expectations
with respect to $\cF_t$ we get that 
\begin{align*}
  \norm{\nabla P(t)}_{H}^2
  +\E\left[\int_{Q_t^T}|\nabla\tilde P|^2
  +\int_t^T\norm{\nabla Z(s)}^2_{\cL^2(K,H)}\,\d s\,\Bigg|\cF_t\right]
  \leq c\E\left[\int_t^T\norm{\nabla P(s)}_H^2\,\d s\,\Bigg|\cF_t\right]\,,
\end{align*}
so that applying the backward stochastic Gronwall 
Lemma~\ref{lem:gronwall} and then taking
expectations yield 
$\nabla \tilde P=0$ almost everywhere in $\Omega\times Q$, 
hence also $\tilde P=0$ almost everywhere in $\Omega\times Q$ since $\tilde P_\OO=0$.
Consequently, the stochastic integral appearing in the estimate above vanishes, 
and we deduce also $\nabla P=0$ in $L^2_\cP(\Omega; C^0([0,T]; H^d))$, 
from which $\tilde P=0$ in $L^2_\cP(\Omega; C^0([0,T]; V_1^*))$.
Also, $\nabla Z=0$ in $L^2_\cP(\Omega; L^2(0,T; \cL^2(K,H^d)))$.
This concludes the proof of Theorem~\ref{thm4}.

%%%%%%%%%%%%%%%%%%%%%%%%%%%%%%%%%%

\section{Necessary conditions for optimality}
\label{sec:nec}
In this last section, we prove the two versions of necessary conditions 
for optimality contained in Theorems~\ref{thm5}--\ref{thm6}.
Let then $\bu\in\mU_{ad}$ be an optimal control for problem {\bf (CP)}
and let us set $(\varphi,\mu):=S(\bu)$ as its corresponding optimal state.
Let us also fix an arbitrary $\bv\in\mU_{ad}$.

By convexity of $\mU_{ad}$ we have $\bu+\delta(\bv-\bu)\in\mU_{ad}$
for all $\delta\in[0,1]$. Hence, setting $(\varphi_\delta, \mu_\delta):=S(\bu+\delta(\bv-\bu))$,
for every $\delta\in[0,1]$ the minimality of $\bu$ yields
\[
  J(\varphi,\bu)\leq
  \frac{\alpha_1}2\E\int_Q|\varphi_\delta - \varphi_Q|^2 
  +\frac{\alpha_2}2\E\int_\OO|\varphi_\delta(T)-\varphi_T|^2 
  +\frac{\alpha_3}{2}\E\int_Q|\bu+\delta (\bv-\bu)|^2\,,
\]
which entails in turn
\begin{align*}
&\frac{\alpha_1}2\E\int_Q\left(|\varphi_\delta|^2 - |\varphi|^2 - 
2(\varphi_\delta-\varphi)\varphi_Q\right)
+\frac{\alpha_2}2\E\int_\OO\left(|\varphi_\delta(T)|^2 
- |\varphi(T)|^2 - 2(\varphi_\delta-\varphi)(T)\varphi_T\right)\\
&\qquad+\frac{\alpha_3}2\E\int_Q\left(\delta^2|\bv-\bu|^2 + 2\delta \bu \cdot(\bv-\bu)\right)
\geq0\,.
\end{align*}
Now, the functions $\zeta\mapsto \E\int_Q|\zeta|^2$ and 
$\zeta\mapsto \E\int_\OO|\zeta|^2$
are Fr\'echet-differentiable on $L^2_\cP(\Omega; L^2(0,T; H))$
and $L^2(\Omega, \cF_T; H)$, respectively.
Hence, the mean-value theorem yields
\begin{align*}
  &\alpha_1\E\int_Q\frac{\varphi_\delta-\varphi}{\delta}\int_0^1
  \left((\varphi + \tau(\varphi_\delta-\varphi))-\varphi_Q\right)\,\d\tau
  +\alpha_3\E\int_Q\bu\cdot(\bv-\bu) + 
  \frac{\alpha_3}2\delta\E\norm{\bv-\bu}^2_{L^2(0,T; H^d)}\\
  &\qquad+\alpha_2\E\int_\OO\frac{\varphi_\delta-\varphi}{\delta}(T)\int_0^1
  \left((\varphi(T) + \tau(\varphi_\delta-\varphi)(T))-\varphi_T\right)\,\d\tau\geq0\,.
\end{align*}
At this point, as $\delta\to0$ we have $\bu+\delta\bv\to\bu$ in $\mU$,
so \eqref{cont}--\eqref{cont2} imply that 
\begin{align*}
  \int_0^1\left((\varphi + \tau(\varphi_\delta-\varphi))-\varphi_Q\right)\,\d\tau \to 
  \varphi - x_Q \qquad&\text{in }
  L^p_\cP(\Omega;  L^2(0,T; V_1))\,,\\
  \int_0^1\left((\varphi(T) + \tau(\varphi_\delta-\varphi)(T))-\varphi_T\right)\,\d\tau \to 
  \varphi(T) - \varphi_T \qquad&\text{in }
  L^{p/3}(\Omega, \cF_T; H)\,.
\end{align*}
Moreover, Theorem~\ref{thm3} ensures that 
\begin{align*}
  \frac{\varphi_\delta-\varphi}{\delta} \wto \theta_{\bv-\bu} \qquad&\text{in }
  L^p_\cP(\Omega; L^2(0,T; H))\,,\\
  \frac{\varphi_\delta-\varphi}{\delta}(T) \wto \theta_{\bv-\bu}(T) \qquad&\text{in }
  L^{p/3}(\Omega, \cF_T; H)\,.
\end{align*}
Hence, noting that $\frac{p}3\geq2$, letting $\delta\to0$ we obtain exactly \eqref{NC},
and Theorem~\ref{thm5} is proved.

Lastly, we note that \eqref{NC2} follows directly from \eqref{NC} provided to show the
duality relation
\[
\alpha_1\E\int_Q\theta_{\bv-\bu}(\varphi-\varphi_Q)+
  \alpha_2\E\int_\OO\theta_{\bv-\bu}(T)(\varphi(T)-\varphi_T)=
  \E\int_Q\varphi(\bv-\bu)\cdot\nabla P\,.
\]
In order to prove this, we can take $g=0$ and $\bh=\bv-\bu$
in the duality relation \eqref{duality_lam}, and then let 
$\lambda\searrow0$ thanks to the convergences 
\eqref{conv1_lin'}--\eqref{conv2_lin'}.
This concludes the proof of Theorem~\ref{thm6}.

\section*{Acknowledgments}
The author is funded by the Austrian Science Fund (FWF)
through project M 2876.
The author also thanks the anonymous referees 
for the positive comments and insightful suggestions that improved
the presentation of the paper.

%%%%%%%%%%%%%%%%%%%%%%%%%%%%%%%%%%%%%%%%%%%

%\bibliography{ref}{}
\bibliographystyle{abbrv}

\def\cprime{$'$}

\end{document}